\documentclass[11pt]{amsart}
\usepackage[T1]{fontenc}
\usepackage[latin1]{inputenc}
\usepackage{a4}
\usepackage{setspace}
\usepackage{xypic}
\usepackage{amssymb}
\usepackage{amsthm}
\usepackage[english]{babel}
\usepackage{latexsym}
\date{}
\newtheorem{thm}{Th\'eor\`eme}[section]
\newtheorem{hypothese}[thm]{Hypothèse}

\newtheorem{defn}[thm]{D\'efinition}
\newtheorem{rem}[thm]{Remarque}
\newtheorem{prop}[thm]{Proposition}
\newtheorem{lem}[thm]{Lemme}
\newtheorem{cor}[thm]{Corollaire}

\newenvironment{f-proof}[1][\sc D\'emonstration.]{\begin{trivlist}
\item[\hskip \labelsep {\bfseries #1}]}{\hfill{$\square$}\end{trivlist}}

\newcommand{\fonc}[5]{
 \begin{array}{cccc}
 #1: & #2 & \longrightarrow & #3\\
     & #4 & \longmapsto & #5
 \end{array}
}
\newcommand{\appl}[4]{
 \begin{array}{cccc}
  #1 & \longrightarrow & #2\\
  #3 & \longmapsto & #4
 \end{array}
}
\begin{document}

\title[Le système dynamique de Deninger et le topos Weil-étale]{Sur l'analogie entre le système dynamique de Deninger et le topos Weil-étale}
\author{Baptiste Morin}

\maketitle

\begin{abstract}
We express some basic properties of Deninger's conjectural dynamical system in terms of morphisms of topoi. Then we show that the current definition of the Weil-étale topos satisfies these properties. In particular, the flow, the closed orbits, the fixed points of the flow and the foliation in characteristic $p$ are well defined on the Weil-étale topos. This analogy extends to arithmetic schemes. Over a prime number $p$ and over the archimedean place of $\mathbb{Q}$, we define a morphism from a topos associated to Deninger's dynamical system to the Weil-étale topos. This morphism is compatible with the structure mentioned above.
\end{abstract}


\section*{Introduction}
C. Deninger a développé un formalisme cohomologique conjectural permettant une
interprétation cohomologique des fonctions $L$ motiviques et des fonctions zêta de schémas arithmétiques. L'existence d'une telle cohomologie permettrait de prouver les grandes conjectures relatives à ces fonctions $L$. Deninger a ensuite montré que certains systèmes dynamiques munis de feuilletages possèdent une cohomologie analogue. Il a finalement suggéré l'existence d'un certain système dynamique feuilleté fonctoriellement attaché à un schéma arithmétique, donnant lieu à la cohomologie de type géométrique attendue. Les propriétés de ce système dynamique ont été étudiées dans une série d'articles (cf.
\cite{Deninger-some-analogies},
\cite{Deninger-possible-significance}, \cite{Deninger-NTandDSonFS},
\cite{Deninger2}, \cite{Deninger-explicit formulas}). La
construction de tels systèmes dynamiques fonctoriellement attachés
aux anneaux d'entiers de corps de nombres par exemple, fait toujours
défaut.

D'autre part, Lichtenbaum a conjecturé dans \cite{Lichtenbaum} l'existence d'une cohomologie, issue d'un certain topos Weil-étale, permettant d'exprimer les valeurs spéciales des fonctions $L$ motiviques et des fonctions zêta de schémas arithmétiques. Une telle cohomologie sera dite de type arithmétique. Nous disposons actuellement d'une définition provisoire du topos Weil-étale d'un schéma arithmétique (cf. \cite{Lichtenbaum} et \cite{Flach-moi}). En suivant le point de vue de \cite{SGA4}, un topos doit être pensé comme un espace topologique généralisé. Le but de cet article est de montrer que certaines propriétés de base du système dynamique de Deninger (flot, orbites fermées, points fixes et feuilletage en caractéristique $p$) sont satisfaites par le topos Weil-étale tel qu'il est défini dans \cite{Flach-moi}.

Dans la deuxième et la troisième section de cet article, nous étudions le
lien existant entre le gros topos Weil-étale $Y_W$ et
le système dynamique conjecturalement associé à une variété $Y$ sur un
corps fini $k=\mathbb{F}_q$. Il est pour cela nécessaire de traduire
les propriétés de ce système dynamique en termes de topos. On montre que le topos Weil-étale est muni d'une projection
$$Y_W\rightarrow B_{W_k}\simeq\mathcal{S}(\mathbb{R},\mathbb{R}/log(q)\mathbb{Z})$$
sur le topos associé à l'espace homogène
$\mathbb{R}/log(q)\mathbb{Z}$. Cette projection est
$\mathbb{R}$-équivariante, en ce sens qu'elle est définie au-dessus
de $B_{\mathbb{R}}$. Une feuille, définie comme la fibre du
morphisme $p$ au-dessus d'un point du cercle $u\in
\mathbb{R}/log(q)\mathbb{Z}$, est équivalente au topos étale
(légèrement modifié) associé au schéma
$\overline{Y}:=Y\times_k\overline{k}$. Un point fermé $v$ de $Y$
fournit une inclusion fermée, à nouveau au-dessus de
$B_{\mathbb{R}}$, du topos associé à
$\mathbb{R}/log(N(v))\mathbb{Z}$ dans le topos Weil-étale. En
supposant que ce système dynamique puisse être fonctoriellement attaché
à $Y$, nous montrons que l'on peut lui associer un topos muni d'un
morphisme canonique dans le topos Weil-étale. Ce morphisme est alors
compatible à l'action de $\mathbb{R}$, aux orbites fermées ainsi
qu'au feuilletage (cf. Théorème \ref{propmorphisme}).

Dans la quatrième et la cinquième section, nous considérons l'anneau d'entiers $\mathcal{O}_K$ d'un corps de nombres, le topos Weil-étale $\bar{X}_W$ (défini dans \cite{Lichtenbaum} et \cite{Flach-moi}) et le système dynamique conjecturalement associé à $\bar{X}=\overline{Spec(\mathcal{O}_K)}$. Ce dernier devrait être un espace laminé
de dimension trois, muni d'un feuilletage, lui-même compatible à
l'action du groupe $\mathbb{R}$. Une place finie du corps de nombres $K$ de
norme $N(v)$ devrait correspondre à une orbite fermée du flot
de longueur $log(N(v))$, alors qu'une place archimédienne
correspondrait à un point fixe. On exprime le flot, les orbites fermées et points fixes de ce système dynamique en termes de morphismes de topos (cf. Théorème \ref{orbiteferme+ptsfixes}), que l'on observe ensuite sur le topos Weil-étale $\bar{X}_W$. Cette étude permet d'interpréter le morphisme canonique $$\mathfrak{f}:\bar{X}_W\rightarrow B_{\mathbb{R}}$$ comme un flot (cf. Section \ref{appendix}). Si $v$ est une place finie de $K$ de norme $N(v)$ (respectivement une place archimédienne), alors le topos Weil-étale $\bar{X}_W$ possède une orbite fermée de longueur $log(N(v))$ (respectivement un point fixe).

Dans la sixième section, nous étendons cette analogie aux schémas arithmétiques de dimension supérieure. Soit $\mathcal{X}$ un schéma connexe, régulier, plat et propre sur $Spec(\mathbb{Z})$. On définit un morphisme de la fibre du système dynamique conjecturalement associé à $\mathcal{X}$ dans la fibre du topos Weil-étale $\bar{\mathcal{X}}_W$ défini dans \cite{Flach-moi}, au dessus d'un point fermé $p\in Spec(\mathbb{Z})$ et au-dessus de la place archimédienne $\infty$ de $\mathbb{Q}$.

Le topos que l'on associe dans ce travail au système dynamique de Deninger est particulièrement simple et maniable. En contrepartie, il réduit ce système dynamique feuilleté à un espace topologique muni d'une action continue de $\mathbb{R}$, et il n'est pas nécessairement connexe sur $\mathcal{T}$. Malgré ses défauts, cette définition retient suffisamment d'information pour observer le flot, les orbites fermées, les points fixes et le feuilletage en caractéristique $p$.\\

{\bf {Remerciements.}} Je tiens à remercier Christopher Deninger, Matthias Flach, Luc Illusie, Masanori Morishita et Frédéric Paugam pour leurs commentaires. Je suis aussi très reconnaissant envers Masanori Morishita pour son hospitalité.

\section{Préliminaires}

\subsection{Topos} Nous rappelons dans cette section quelques notions simples (issues de \cite{SGA4} IV) concernant les topos et leurs morphismes, que nous utiliserons tout au long de ce travail.

Soit $C$ une catégorie munie d'une topologie de Grothendieck $J$. On note $\widetilde{(C,J)}$ la catégorie des faisceaux d'ensembles sur le site $(C,J)$ (nous passons les questions d'univers sous silence). Un topos est une catégorie $\mathcal{S}$ équivalente à une catégorie de la forme $\widetilde{(C,J)}$. Les axiomes de Giraud donnent une caractérisation des topos (cf. \cite{SGA4} IV.1). Un morphisme de topos $f:\mathcal{S}'\rightarrow \mathcal{S}$ est la donnée d'un couple de foncteurs adjoints $f^*:\mathcal{S}\rightarrow \mathcal{S}'$ et $f_*:\mathcal{S}'\rightarrow \mathcal{S}$, de sorte que l'adjoint à gauche $f^*$ de $f_*$ soit exact à gauche. Rappelons qu'un foncteur est exact à gauche s'il commute aux limites projectives finies (i.e. s'il préserve l'objet final et s'il commute aux produits fibrés).

Un site $(C,J)$ est dit \emph{exact à gauche} si les produits fibrés et un objet final existent dans $C$, et si $J$ est moins fine que la topologie canonique. Cette dernière condition signifie simplement que les préfaisceaux représentables $y(X):=Hom_C(-,X)$ sont des faisceaux. Si $(C,J)$ est un site exact à gauche alors le foncteur
$$y:C\rightarrow \widetilde{(C,J)},$$
envoyant $X$ sur $y(X)$, est pleinement fidèle. C'est le \emph{plongement de Yoneda}. Ce foncteur commute par définition aux limites projectives quelconques. En particulier, l'objet final de $\widetilde{(C,J)}$ est le faisceau représenté par l'objet final de $C$.

Soient  $(C,J)$ et  $(C',J')$ deux sites exacts à gauche. Un \emph{morphisme de sites exacts à gauche}
$f^*:(C,J)\rightarrow (C',J')$
est un foncteur $f^*$ exact à gauche qui est de plus continu (une famille couvrante dans $C$ est envoyée sur une famille couvrante dans $C'$). Un tel morphisme de sites  exacts à gauche induit un morphisme de topos  (cf. \cite{SGA4} IV.4.9.2)
$$f:=(f^*,f_*):\widetilde{(C',J')}\longrightarrow\widetilde{(C,J)}$$
Soient $\mathcal{S}$ et $\mathcal{S}'$ deux topos et soit $f^*:\mathcal{S}\rightarrow\mathcal{S}'$ un foncteur commutant aux limites projectives finies et aux limites inductives quelconques. Alors il existe  (essentiellement) un unique morphisme de topos $f:\mathcal{S}'\rightarrow\mathcal{S}$ dont l'image inverse est $f^*$. Par exemple, soit $\mathcal{S}$ un topos et $X$ un objet de $\mathcal{S}$. Alors $\mathcal{S}/X$ est un topos, que l'on appelle \emph{topos induit}. Le foncteur changement de base
$$
\appl{\mathcal{S}}{\mathcal{S}/X}{F}{F\times X}
$$
commute aux limites projectives et inductives quelconques (car ces limites sont universelles dans $\mathcal{S}$). On obtient un morphisme
$$\mathcal{S}/X\rightarrow \mathcal{S}$$
que l'on appelle \emph{morphisme de localisation}.

Un morphisme de topos $i:\mathcal{S}'\rightarrow\mathcal{S}$ est un \emph{plongement} si l'image directe $i_*$ est pleinement fidèle. Une sous-catégorie $\mathcal{E}$ de $\mathcal{S}$ est un \emph{sous-topos} lorsqu'il existe un plongement $i:\mathcal{S}'\rightarrow\mathcal{S}$ tel que $\mathcal{E}$ soit l'image essentielle de $i_*$. En d'autre termes, un sous-topos  $\mathcal{E}$ de $\mathcal{S}$ est une sous-catégorie strictement pleine telle que le foncteur d'inclusion $\mathcal{E}\rightarrow\mathcal{S}$ est l'image directe d'un morphisme de topos, i.e. telle que le foncteur d'inclusion admette un adjoint à gauche qui est exact à gauche. Soit $f:\mathcal{E}\rightarrow\mathcal{S}$ un morphisme. Alors il existe un unique sous-topos $Im(f)$ de $\mathcal{S}$ tel qu'il y ait une factorisation $\mathcal{E}\rightarrow Im(f)\rightarrow\mathcal{S}$ où $\tilde{f}:\mathcal{E}\rightarrow Im(f)$ est \emph{surjectif}, i.e. $\tilde{f}^*$ est fidèle. Le sous-topos $Im(f)$ de $\mathcal{S}$ est \emph{l'image du morphisme $f$}.

Soit $\mathcal{S}$ un topos et $U$ un sous-objet de l'objet final de $\mathcal{S}$. Alors on dit que le morphisme $j:\mathcal{S}/U\rightarrow \mathcal{S}$ est un \emph{plongement ouvert}. Le \emph{sous-topos fermé complémentaire} est la sous-catégorie strictement pleine $F$ de $\mathcal{S}$ constituée des objets $X$ de $\mathcal{S}$ tels que $j^*X$ est l'objet final de $\mathcal{S}/U$. Un \emph{plongement fermé} est un morphisme
$i:\mathcal{E}\rightarrow\mathcal{S}$
tel que $i_*$ induise une équivalence de $\mathcal{E}$ sur un sous-topos fermé de $\mathcal{S}$.

Nous dirons qu'un diagramme de topos
\[ \xymatrix{
\mathcal{E}'\ar[d]_{i'}\ar[r]&\mathcal{E} \ar[d]_{i}   \\
\mathcal{S}'\ar[r]^f& \mathcal{S}
} \] est un \emph{pull-back} lorsqu'il est commutatif et 2-cartésien (cf. \cite{SGA4} IV. Proposition 5.11). Un tel diagramme commutatif est 2-cartésien précisément lorsque le morphisme $\mathcal{E}'\rightarrow\mathcal{S}'\times_{\mathcal{S}}\mathcal{E}$ est une équivalence, où $\mathcal{S}'\times_{\mathcal{S}}\mathcal{E}$ est un 2-produit fibré (les 2-produits fibrés existent dans la 2-catégorie des topos). Considérons un pull-back comme ci-dessus tel que $i:\mathcal{E}\rightarrow\mathcal{S}$ soit un plongement fermé (resp. ouvert). Alors $i':\mathcal{E}'\rightarrow\mathcal{S}'$ est un plongement fermé (resp. ouvert). En particulier, l'image inverse d'un sous-topos fermé (resp. ouvert) est un sous-topos fermé (resp. ouvert). Le sous-topos $Im(i')$ est \emph{l'image inverse} du sous-topos  $Im(i)$ par le morphisme $f:\mathcal{S}'\rightarrow\mathcal{S}$.

Le lecteur peu familier avec la théorie des topos pourra voir un topos simplement comme un espace topologique généralisé, un morphisme de topos comme une application continue, un morphisme de localisation comme un homéomorphisme local, un plongement comme l'immersion d'un sous-espace topologique, un plongement ouvert (resp. fermé) de topos comme un plongement plongement ouvert (resp. fermé) d'espaces et un 2-produit fibré de topos comme un produit fibré d'espaces topologiques.

\subsection{Le gros topos d'un espace topologique}\label{subsect-local-section}

Soit $({Top},\mathcal{J}_{ouv})$ le site constitué de la
catégorie des espaces topologiques,
munie de la topologie engendrée par la prétopologie des
recouvrements ouverts : $\{X_i\rightarrow X, i\in I\}\in Cov(X)$ si $X=\cup_{i\in I} X_i$ est un recouvrement ouvert. On note $\mathcal{J}_{ls}$ la topologie sur
${Top}$ engendrée par la prétopologie des recouvrements
admettant des sections locales : $\{X_i\rightarrow X, i\in I\}\in Cov(X)$ si pour tout $x\in X$, il existe $i\in I$, un voisinage ouvert $x\in U\subset X$ et une section continue $U\rightarrow X_i$ de la flèche $X_i\rightarrow X$ au-dessus de $U$. Ces deux topologies sont en fait les mêmes
(i.e. $\mathcal{J}_{ls}=\mathcal{J}_{ouv}$) et nous désignons par
$\mathcal{T}$ le topos des faisceaux sur le site $({Top},\mathcal{J}_{ouv})$. Soit $X$ un espace topologique. Il représente un faisceau $y(X)=Hom_{Top}(-,X)$ sur le site $({Top},\mathcal{J}_{ouv})$. On considère alors la catégorie $\mathcal{T}/_{y(X)}$ des objets de $\mathcal{T}$ au-dessus de $y(X)$. On définit la topologie $\mathcal{J}_{ouv}$ des recouvrements ouverts sur la catégorie $Top/_{X}$ des espaces topologiques au-dessus de $X$. Alors on a une équivalence
$$\mathcal{T}/_{y(X)}\simeq\widetilde{({Top}/_{X},\mathcal{J}_{ouv})}$$
où $\widetilde{({Top}/_{X},\mathcal{J}_{ouv})}$ désigne la catégorie des faisceaux d'ensembles sur le site $({Top}/_{X},\mathcal{J}_{ouv})$. Le topos $TOP(X):=\widetilde{({Top}/_{X},\mathcal{J}_{ouv})}$ est le \emph{gros topos de l'espace $X$} (cf. \cite{SGA4} IV.2.5). En particulier, $\mathcal{T}/_{y(X)}$ est un topos, équivalent au gros topos de l'espace $X$.
\begin{defn}
On note
$$\mathcal{S}(X):=\mathcal{T}/_{y(X)}\simeq TOP(X)$$
le gros topos de l'espace $X$.
\end{defn}
On note $Sh(X)$ la catégorie des espaces étalés sur $X$. On a alors un morphisme $$\mathcal{S}(X)\longrightarrow Sh(X),$$ qui possède une section, et ces deux topos sont cohomologiquement équivalents (cf. \cite{SGA4} IV.2.5).

\subsection{Le gros topos d'un espace muni de l'action d'un groupe topologique} Soit $G$ un groupe topologique.
Le plongement de Yoneda étant exact à gauche, le groupe topologique $G$ représente un groupe
$y(G)$ de $\mathcal{T}$, i.e. un faisceau de groupes sur $({Top},\mathcal{J}_{ouv})$.
\begin{defn}\label{defn-topos-class}
Le \emph{topos classifiant}
$B_{G}$ du groupe topologique $G$ est défini comme
la catégorie des objets de $\mathcal{T}$ munis d'une action à gauche
de $y(G)$.
\end{defn}
Un objet de $B_{G}$ est donc un faisceau d'ensembles $\mathcal{F}$ sur $({Top},\mathcal{J}_{ouv})$, tel que pour tout espace $T$ de $Top$, on ait une action
$$Hom_{Top}(T,G)\times \mathcal{F}(T)\rightarrow\mathcal{F}(T)$$
fonctorielle en $T$. Soit ${Top}^{G}$ la catégorie
des espaces topologiques sur lesquels $G$ opère
continûment. La topologie (encore notée $\mathcal{J}_{ls}$) des
recouvrements admettant des sections locales est définie comme la
topologie induite sur ${Top}^{G}$ par le foncteur
d'oubli ${Top}^{G}\rightarrow {Top}$. Le
plongement de Yoneda définit un foncteur pleinement fidèle
$${Top}^{G}\longrightarrow B_{G}.$$
Via ce foncteur, $\mathcal{J}_{ls}$ est la topologie induite par la
topologie canonique de $B_{G}$ et
${Top}^{G}$ est une sous-catégorie génératrice de
$B_{G}$. On en déduit une équivalence de topos
$$B_{G}\longrightarrow\widetilde{({Top}^{G};\mathcal{J}_{ls})}.$$
Soit maintenant $(G,Z)$ une action continue du groupe
topologique $G$ sur un espace $Z$. L'objet
$y(G,Z)$, défini par le plongement de Yoneda, est un objet
de $B_{G}$. On peut donc définir le topos induit
$\mathcal{S}(G,Z):=B_{G}/_{y(G,Z)}$.
D'après (\cite{SGA4} III Proposition 5.4), on a une équivalence
$$\mathcal{S}(G,Z):=B_{G}/_{y(G,Z)}\simeq\widetilde{({Top}^{G}/_{(G,Z)};\mathcal{J})}$$
où $\mathcal{J}$ est la topologie induite par $\mathcal{J}_{ls}$ sur
${Top}^{G}$ via le foncteur
$${Top}^{G}/_{(G,Z)}\rightarrow{Top}^{G},$$
qui consiste à oublier la flèche sur $(G,Z)$. Cette
topologie $\mathcal{J}$ est donc induite par la topologie des
sections locales sur ${Top}$. On la note à nouveau
$\mathcal{J}_{ls}$.
\begin{prop}\label{prop-local-section-site-general}
Le site
$({Top}^{G}/_{(G,Z)},\mathcal{J}_{ls})$
est un site pour le topos $\mathcal{S}(G,Z):=B_{G}/_{y(G,Z)}$.
\end{prop}
Dans la suite, on appelle $\mathcal{S}(G;Z):=B_{G}/_{y(G;Z)}$
\emph{le topos des gros $G$-faisceaux sur $Z$}. Cette
terminologie se justifie par le fait qu'un objet de
$\mathcal{S}(G;Z)$ est donné par un objet $\mathcal{F}$ de
$\mathcal{T}/_{y(Z)}$ muni d'une action de $y(G)$ telle que
le diagramme
\[ \xymatrix{
 y(G)\times\mathcal{F} \ar[d]_{} \ar[r]^{} &\mathcal{F} \ar[d]_{}   \\
 y(G)\times y(Z)  \ar[r]^{}& y(Z)
} \] soit commutatif. L'équivalence
$$\mathcal{T}/_{y(Z)}:=\widetilde{({Top},\mathcal{J}_{ls})}/_{y(Z)}
\simeq\widetilde{({Top}/_{Z},\mathcal{J}_{ls})}=
TOP(Z),$$ où $TOP(Z)$ désigne le gros topos de l'espace topologique
$Z$ (cf. \cite{SGA4} IV.2.5), montre alors que $\mathcal{F}$ est un
gros faisceau sur $Z$ muni d'une action de $G$ compatible à
celle définie sur $Z$.
\begin{defn}
Le \emph{topos des gros $G$-faisceaux sur $Z$} est le topos induit
$$\mathcal{S}(G;Z):=B_{G}/_{y(G;Z)}\simeq\widetilde{({Top}^{G}/_{(G,Z)},\mathcal{J}_{ls})}$$
\end{defn}
Cette définition sera utilisée en particulier dans la situation suivante. Lorsqu'un espace $Z$ est muni d'une action continue de du groupe $\mathbb{R}$, on note
$$\mathcal{S}(\mathbb{R};Z):=B_{\mathbb{R}}/_{y(\mathbb{R};Z)}\simeq\widetilde{({Top}^{\mathbb{R}}/_{(\mathbb{R},Z)},\mathcal{J}_{ls})}.$$
On utilisera aussi le topos suivant.
\begin{defn}
Soit $G$ un groupe discret ou profini. Le petit topos classifiant $B_G^{sm}$ est la catégorie des ensembles sur lesquels le groupe $G$ opère continûment.
\end{defn}

\section{Le système dynamique en caractéristique positive}\label{sect-systdyn-charp}

Soit $Y$ un schéma lisse, séparé et de type fini sur un corps
fini $k=\mathbb{F}_q$. On note $d=dim(Y)$. La donnée de $Y$ est équivalente à celle du
couple $(\overline{Y}=Y\otimes\overline{k};\varphi)$, où $\varphi^{-1}=F\otimes Id_{\overline{k}}$
avec $F$ le Frobenius absolu de ${Y}$. Alors
l'ensemble des points fermés $Y^0$ de $Y$ est en bijection avec
celui des orbites finies $\mathfrak{o}$ de l'ensemble
$\overline{Y}(\overline{k})$, des $\overline{k}$-points de
$\overline{Y}$ sous l'action de $\varphi$. Si le point fermé $y\in
Y^0$ correspond à l'orbite $\mathfrak{o}$, alors
$$log (N(y))=|\mathfrak{o}|\,log (q).$$ Le couple $(\overline{Y};\varphi)$
peut être vu comme l'analogue d'un couple $(\mathrm{N},\varphi)$,
où $\mathrm{N}$ est un espace  de dimension
$2d$, dont $\varphi$ est un automorphisme.
Plus précisément, $\mathrm{N}$ est un espace laminé compact de la forme décrite dans \cite{Leichtnam-invitation} Definition 2. En particulier, il existe un espace totalement discontinu $\Omega$ tel que $\mathrm{N}$ soit localement homéomorphe à un produit $U\times T$, où $D\subset\mathbb{C}^d$ est un ouvert non-vide et $T\subset \Omega$ est ouvert. $\mathrm{N}$ est donc de dimension topologique $2d$. On note $W_k$ le groupe $\{\varphi^{\mathbb{Z}}\}$ et
$$
\appl{W_k}{\mathbb{R}^*_+}{\varphi}{q}
$$
le morphisme canonique. En suivant C. Deninger \cite{Deninger3} et E. Leichtnam \cite{Leichtnam-invitation} Hypothesis 1, l'action de $W_k$ sur $\mathrm{N}$ permet de
définir une action de $\mathbb{R}$ sur un espace $\mathrm{M}$ de dimension $2\,d+1$. En effet, $W_k$ opère sur le
produit $\mathbb{R}^*_+\times\mathrm{N}$ par la formule (cf. \cite{SGA4} IV.4.5.1)
\begin{eqnarray*}
\varphi^{\nu}\,.\,(u;n)&=&(u.q^{-\nu};\varphi^{\nu}(n))
\end{eqnarray*}
pour $\nu\in\mathbb{Z}$, $u\in\mathbb{R}_+^*$ et $n\in\mathrm{N}$. Autrement dit, l'action de $W_k$ se fait à droite sur $\mathrm{N}$
et à gauche sur $\mathbb{R}$. On définit alors le quotient
$$\mathrm{M}:=(\mathbb{R}^*_+\times\mathrm{N})/W_k.$$
Le groupe $\mathbb{R}$ opère sur $\mathrm{M}$ par
la formule $\phi^t[u;n]=[e^tu;n]$. Dans cette situation, les orbites fermées $\gamma$ de l'action de $\mathbb{R}$
sur $\mathrm{M}$ correspondent bijectivement aux orbites finies
$\mathfrak{o}$ de $\mathrm{N}$ sous l'action de
$W_k=\{\varphi^{\mathbb{Z}}\}$. La longueur d'une telle orbite est
donnée par $l(\gamma)=|\mathfrak{o}|\,log (q)$. D'autre part, la projection canonique
$$p:\mathrm{M}=(\mathbb{R}^*_+\times\mathrm{N})/W_k\longrightarrow
\mathbb{R}^*_+/W_{k}$$ est $\mathbb{R}$-équivariante,
lorsque $\mathbb{R}\simeq \mathbb{R}^*_+$ opère naturellement sur l'espace homogène
$\mathbb{R}^*_+/W_{k}$. Les fibres
$\mathcal{F}_{\overline{u}}:=p^{-1}(\{\overline{u}\})$ de $p$
définissent un feuilletage $\mathcal{F}$ de codimension $1$ sur
$\mathrm{M}$ compatible à l'action de $\mathbb{R}$. Explicitement,
les feuilles de $\mathcal{F}$ sont les images des immersions fermées
$$
\appl{\mathrm{N}}{\mathrm{M}}{n}{[u;n]}
$$
pour $u\in\mathbb{R}^*_+$. Les trajectoires du flot sont partout
perpendiculaires aux feuilles. Dans cette section, nous traduisons
la situation décrite ci-dessus en termes de topos. Nous adoptons une notation additive (via le logarithme $\mathbb{R}^*_+\simeq \mathbb{R}$).

\subsection{Le système dynamique associé à un corps fini}

Soient $\mathbb{F}_q$ un corps fini et
$\overline{\mathbb{F}}_q/\mathbb{F}_q$ une clôture algébrique. En suivant Deninger, on
associe à $Spec(\mathbb{F}_q)$ le système dynamique
$\mathbb{R}/log(q)\mathbb{R}$ où $\mathbb{R}$ opère par translations
à gauche (cf. \cite{Deninger-possible-significance} 2.7).

\begin{thm}\label{morphisme alpha}
Le foncteur $\alpha^*$ qui associe le système dynamique
$\mathbb{R}/log(q)\mathbb{Z}$ au corps fini $\mathbb{F}_q$ induit le morphisme de
topos
$$\alpha:B_{W_{\mathbb{F}_q}}\longrightarrow
B^{sm}_{G_{\mathbb{F}_q}}$$ donné par la flèche
canonique $W_{\mathbb{F}_q}\rightarrow G_{\mathbb{F}_q}$.
\end{thm}

\begin{f-proof}
On note $$\mathbb{M}_q:=(\mathbb{R},\mathbb{R}/log(q)\mathbb{Z})$$ le système
dynamique associé à $\mathbb{F}_q$. La clôture algébrique $\overline{\mathbb{F}}_q/\mathbb{F}_q$ correspond à un point marqué sur $\mathbb{M}_q$. Une extension finie $\mathbb{F}_{q^n}/\mathbb{F}_q$ de degré $n$
induit un morphisme $\mathbb{R}$-équivariant
$$\mathbb{R}/log(q^n)\mathbb{Z}\longrightarrow\mathbb{R}/log(q)\mathbb{Z}.$$
Cette flèche est d'ailleurs un revêtement étale de degré $n$.  On note
$W_{\mathbb{F}_q}$ l'image dans $\mathbb{R}$ du groupe de Weil de
$\mathbb{F}_q$. Alors $\mathbb{M}_q$ est l'espace homogène
$(\mathbb{R},\mathbb{R}/W_{\mathbb{F}_q})$, avec
$W_{\mathbb{F}_q}=log(q)\mathbb{Z}\subset\mathbb{R}.$ On considère le foncteur suivant :
\begin{equation}\label{foncteur-alpha}
Et_{Spec(\mathbb{F}_q)}\simeq{G_{\mathbb{F}_q}}-\underline{Set}^f\longrightarrow W_{\mathbb{F}_q}-\underline{Set}\longrightarrow B_{W_{\mathbb{F}_q}}\simeq B_{\mathbb{R}}/y(\mathbb{M}_q)
\end{equation}
Ci-dessus, l'équivalence entre la catégorie $Et_{Spec(\mathbb{F}_q)}$ des $Spec(\mathbb{F}_q)$-schémas étales et la catégorie des
$G_{\mathbb{F}_q}$-ensembles finis $G_{\mathbb{F}_q}-\underline{Set}^f$ est donnée par la clôture algébrique $\overline{\mathbb{F}}_q/\mathbb{F}_q$. La catégorie $W_{\mathbb{F}_q}-\underline{Set}$ est celle des ensembles munis d'une action du groupe de Weil $W_{\mathbb{F}_q}$. Le foncteur $G_{\mathbb{F}_q}-\underline{Set}^f\rightarrow W_{\mathbb{F}_q}-\underline{Set}$ est induit par le morphisme canonique $W_{\mathbb{F}_q}\rightarrow G_{\mathbb{F}_q}$ Le foncteur $W_{\mathbb{F}_q}-\underline{Set}\rightarrow B_{W_{\mathbb{F}_q}}$ envoie un $W_{\mathbb{F}_q}$-ensemble $E$ sur le faisceau de $B_{W_{\mathbb{F}_q}}$ représenté par le $W_{\mathbb{F}_q}$-espace discret $E$. Enfin, l'équivalence $B_{W_{\mathbb{F}_q}}\simeq B_{\mathbb{R}}/y(\mathbb{M}_q)$ est donnée dans \cite{SGA4} IV.5.8 (voir aussi la section \ref{subsect-equivalence} ci-dessous). L'image essentielle du foncteur (\ref{foncteur-alpha}) est contenue dans la sous-catégorie pleine $$Top^{\mathbb{R}}/_{\mathbb{M}_{q}}\hookrightarrow B_{\mathbb{R}}/y(\mathbb{M}_q).$$
Le foncteur défini ci-dessus est de la forme suivante :
$$
\fonc{\alpha^*}{Et_{Spec(\mathbb{F}_q)}}{Top^{\mathbb{R}}/_{\mathbb{M}_{q}}}
{Spec(\mathbb{F}_{q^{n_1}}\times...\times\mathbb{F}_{q^{n_s}})}
{\mathbb{M}_{q^{n_1}}\coprod...\coprod\mathbb{M}_{q^{n_s}}}
$$
Ce foncteur est exact à gauche, puisqu'il préserve l'objet final
et les produits fibrés. Un recouvrement étale de $Et_{Spec(\mathbb{F}_q)}$ est envoyé sur un
recouvrement de ${Top}^{\mathbb{R}}/_{\mathbb{M}_{q}}$
pour la topologie des sections locales, puisqu'un étalement (i.e. un homéomorphisme local) admet
des sections locales au-dessus de son image. Ainsi, le foncteur
$$\alpha^*:(Et_{Spec(\mathbb{F}_q)};\mathcal{J}_{et})\longrightarrow({Top}^{\mathbb{R}}/_{\mathbb{M}_{q}};\mathcal{J}_{ls})$$
est un morphisme de sites exacts à gauche, puisqu'il est continu et
exact à gauche. On en déduit l'existence d'un morphisme de topos
$$\alpha:\widetilde{({Top}^{\mathbb{R}}/_{\mathbb{M}_{q}};\mathcal{J}_{ls})}\longrightarrow
\widetilde{(Et_{Spec(\mathbb{F}_q)};\mathcal{J}_{et})}.$$ A nouveau d'après \cite{SGA4} III Proposition 5.4 et \cite{SGA4} IV.5.8, on a les équivalences
$$\widetilde{({Top}^{\mathbb{R}}/_{\mathbb{M}_{q}};\mathcal{J}_{ls})}
\simeq
{B_{\mathbb{R}}}/_{y(\mathbb{R},\mathbb{R}/W_{\mathbb{F}_q})}\simeq
B_{W_{\mathbb{F}_q}}.$$ De plus, le topos étale de
$Spec(\mathbb{F}_q)$ s'identifie à la catégorie des ensembles sur
lesquels le groupe de Galois $G_{\mathbb{F}_q}$ opère continûment,
une clôture séparable $\overline{\mathbb{F}_q}/\mathbb{F}_q$ ayant
été choisie. En d'autres termes, on a une équivalence
$$\widetilde{(Et_{Spec(\mathbb{F}_q)};\mathcal{J}_{et})}\simeq B^{sm}_{G_{\mathbb{F}_q}},$$
où $B^{sm}_{G_{\mathbb{F}_q}}$ est le petit topos classifiant du
groupe profini $G_{\mathbb{F}_q}$. On obtient le morphisme
$$B_{W_{\mathbb{F}_q}}\simeq \widetilde{({Top}^{\mathbb{R}}/_{\mathbb{M}_{q}};\mathcal{J}_{ls})}
\longrightarrow\widetilde{(Et_{Spec(\mathbb{F}_q)};\mathcal{J}_{et})}\simeq B^{sm}_{G_{\mathbb{F}_q}}$$
\end{f-proof}

Le résultat précédent établit un lien entre l'idée de
Deninger, consistant à voir un corps fini comme le système dynamique
$\mathbb{M}_q$, et celle de Lichtenbaum qui consiste à remplacer
le groupe de Galois d'un corps fini par son groupe de Weil. Le topos Weil-étale $B_{W_{\mathbb{F}_q}}$ de $Spec(\mathbb{F}_q)$ est précisément le gros topos associé à l'espace homogène
$\mathbb{M}_q=\mathbb{R}/log(q)\mathbb{Z}$.

\subsection{L'équivalence $\mathcal{S}(W_k;\mathrm{N})\simeq\mathcal{S}(\mathbb{R};\mathrm{M})$}\label{subsect-equivalence}
Nous reprenons les notations de l'introduction de cette section \ref{sect-systdyn-charp}.
\begin{prop}
Les topos $\mathcal{S}(W_k;\mathrm{N})$ et
$\mathcal{S}(\mathbb{R};\mathrm{M})$ sont canoniquement équivalents
en tant que topos sur $B_{\mathbb{R}}$.
\end{prop}
\begin{f-proof}
Le morphisme $l:W_k=W_{\mathbb{F}_q}\rightarrow\mathbb{R}$
définit un morphisme de topos classifiants $B_l:B_{W_k}\rightarrow B_{\mathbb{R}}$.
Le pull-back $B_l^*$ est défini par restriction du groupe
d'opérateurs, il commute aux limites inductives et projectives
quelconques. Ce foncteur possède donc un adjoint à droite et un
adjoint à gauche. On obtient une suite de trois foncteurs adjoints
$$B_{l!}\,\,\,;\,\,\,\,\,\,B_l^*\,\,\,;\,\,\,\,\,\,B_{l*}.$$
Le foncteur
$$
\fonc{B_{l!}}{B_{W_k}}{B_{\mathbb{R}}}{Z}{y(\mathbb{R})\times^{y(W_k)}Z}
$$
envoie un objet $Z$ de $B_{W_k}$ sur le quotient
$y(\mathbb{R})\times^{y(W_k)}Z:=(y(\mathbb{R})\times Z)/y(W_k)$, où
$y(W_k)$ opère à gauche sur $Z$ et à droite sur $y(\mathbb{R})$ via
le morphisme $l$. Alors $B_{l!}$ se factorise à travers le foncteur
$$
\fonc{f_!}{B_{\mathbb{R}}/_{y(\mathbb{R},\mathbb{R}/W_k)}}{B_{\mathbb{R}}}{Z\rightarrow
y(\mathbb{R}/W_k)}{Z}.
$$
Le foncteur $f_!$ est adjoint à gauche de l'image
inverse du morphisme de localisation
$f:B_{\mathbb{R}}/_{y(\mathbb{R}/W_k)}\rightarrow
B_{\mathbb{R}}$. Ce qui précède induit l'équivalence
$h:B_{W_k}\simeq B_{\mathbb{R}}/_{y(\mathbb{R},\mathbb{R}/W_k)}$
dont l'image directe est donnée par le foncteur
$$
\fonc{h_*}{B_{W_k}}{B_{\mathbb{R}}/_{y(\mathbb{R},\mathbb{R}/W_k)}}{Z}{y(\mathbb{R}\times^{y(W_k)}Z)\rightarrow
y(\mathbb{R}/W_k)}.
$$
Les morphismes $B_l$ et $f\circ h$ sont canoniquement
isomorphes dans la catégorie
$\underline{Homtop}\,(B_{W_k};B_{\mathbb{R}})$. Ainsi, le morphisme
$B_l$ s'interprète comme le morphisme de localisation
$$B_{W_k}\simeq B_{\mathbb{R}}/_{y(\mathbb{R},\mathbb{R}/W_k)}\longrightarrow
B_{\mathbb{R}}.$$
On conserve les notations de l'introduction de cette section \ref{sect-systdyn-charp}.
Le topos des gros $W_{k}$-faisceaux sur $\mathrm{N}$ est défini comme le topos induit
$\mathcal{S}(W_k;\mathrm{N}):=B_{W_k}/_{y(W_k;\mathrm{N})}$, où l'objet $y(W_k;\mathrm{N})$ de $B_{W_k}$ est donné par l'action de
$W_k$ sur $\mathrm{N}$. On a donc un morphisme de localisation
$$\mathcal{S}(W_k;\mathrm{N})\longrightarrow B_{W_k}$$
qui traduit l'action de $W_k$ sur $\mathrm{N}$. L'équivalence $h$
induit une équivalence de topos induits
$$B_{W_k}/_{y(W_k,\mathrm{N})}\longrightarrow (B_{\mathbb{R}}/_{y(\mathbb{R}/W_k)})/_{h_*y(W_k,\mathrm{N})}
=B_{\mathbb{R}}/_{h_*y(W_k,\mathrm{N})}.$$
L'objet $h_*y(W_k,\mathrm{N})$
est représenté par l'action de $\mathbb{R}$ sur
$\mathrm{M}:=(\mathbb{R}\times\mathrm{N})/W_k$. On obtient un
diagramme commutatif de topos
 \[ \xymatrix{
\mathcal{S}(W_k;\mathrm{N})\ar[d]_{}\ar[r]^{}& \mathcal{S}(\mathbb{R};\mathrm{M})\ar[d]^{p}\\
 B_{W_k}\ar[dr]_{B_l}\ar[r]^{h}
 &B_{\mathbb{R}}/_{y(\mathbb{R},\mathbb{R}/W_k)}\ar[d]^{f}\\
& B_{\mathbb{R}} }
\]
où les flèches horizontales sont des équivalences.
\end{f-proof}
\subsection{Le morphisme structural et les orbites fermées}
D'après ce qui précède, on a un morphisme
$p:\mathcal{S}(\mathbb{R};\mathrm{M})\rightarrow B_{\mathbb{R}}/_{y(\mathbb{R},\mathbb{R}/W_k)}$.

\begin{defn}
Le \emph{morphisme flot} $\mathfrak{f}$ est le morphisme de
localisation
$$\mathfrak{f}:\mathcal{S}(\mathbb{R};\mathrm{M})=B_{\mathbb{R}}/_{y(\mathbb{R};\mathrm{M})}\longrightarrow
B_{\mathbb{R}}.$$
\end{defn}
On voit d'ailleurs que le morphisme $\mathfrak{f}$ se factorise à
travers $p$. Soit $\gamma$ une orbite fermée de l'action de
$\mathbb{R}$ sur $\mathrm{M}$ correspondant à une orbite finie
$\mathfrak{o}$ de $\mathrm{N}$ sous l'action de $W_k$. En posant
$|\mathfrak{o}|=n$, on a $l(\gamma)=|\mathfrak{o}|\,log (q)=n\,log (q)$. On pose $k(\gamma)=\mathbb{F}_{q^n}$.

\begin{prop}\label{orbite-ferme-syst-dyn-char-p}
Une orbite fermée $\gamma$ induit un plongement fermé
$$i_{\gamma}:B_{W_{k(\gamma)}}\longrightarrow\mathcal{S}(\mathbb{R};\mathrm{M})$$
tel que la composition
$\mathfrak{f}\circ
i_{\gamma}:B_{W_{k(\gamma)}}\rightarrow B_{\mathbb{R}}$ est
induite par le morphisme canonique
$W_{k(\gamma)}\rightarrow\mathbb{R}$.
\end{prop}
\begin{f-proof}
L'immersion fermée
$\gamma:\mathbb{M}_{q^n}:=\mathbb{R}/nlog(q)\mathbb{Z}\longrightarrow\mathrm{M}$
est $\mathbb{R}$-équivariante. Elle induit un plongement fermé de
topos
$$i_{\gamma}:B_{\mathbb{R}}/_{y(\mathbb{M}_{q^n})}\longrightarrow\mathcal{S}(\mathbb{R};\mathrm{M})$$
tel que la composition
$p\circ
i_{\gamma}:B_{\mathbb{R}}/_{y(\mathbb{M}_{q^n})}\rightarrow
B_{\mathbb{R}}/_{y(\mathbb{M}_{q})}$ soit donnée par le revêtement
$\mathbb{R}$-équivariant
$$\mathbb{M}_{q^n}=\mathbb{R}/(|\mathfrak{o}|\,log
(q)\mathbb{Z})\longrightarrow\mathbb{R}/log
(q)\mathbb{Z}=\mathbb{M}_{q}.$$ On note $k(\gamma)=\mathbb{F}_{q^n}$.
Le morphisme $i_{\gamma}$ est un plongement fermé par le lemme \ref{plongementferme-equi}. On obtient le résultat en utilisant les équivalences
$B_{\mathbb{R}}/_{y(\mathbb{M}_{q^n})}=B_{W_{k(\gamma)}}$ et
$B_{\mathbb{R}}/_{y(\mathbb{M}_{q})}=B_{W_{k}}$ (voir la sous-section
\ref{subsect-orbites-fermee-char-0} pour plus de détails sur ce qui
précède).
\end{f-proof}

\subsection{Le topos des faisceaux sur $\mathrm{M}$ est un produit fibré}
On note $E_{\mathbb{R}}$ l'objet de $B_{\mathbb{R}}$ représenté par
l'espace topologique $\mathbb{R}$ sur lequel $\mathbb{R}$ opère par
translations à gauche. Alors on a un isomorphisme $\mathcal{T}\simeq
B_{\mathbb{R}}/_{E_{\mathbb{R}}}$ de sorte que la composition
$\mathcal{T}\simeq B_{\mathbb{R}}/_{E_{\mathbb{R}}}\rightarrow
B_{\mathbb{R}}$ soit induite par le morphisme de groupes
$\{1\}\rightarrow\mathbb{R}$ (cf \cite{SGA4} IV.5.8.2). On note
aussi
$$\mathcal{S}(\mathrm{M}):=\mathcal{T}/_{y(\mathrm{M})}\simeq\widetilde{({Top}/_{\mathrm{M}},\mathcal{J}_{ouv})}:=TOP(\mathrm{M})$$
le gros topos de l'espace topologique $\mathrm{M}$ (cf \cite{SGA4}
IV.2.5).
\begin{lem}\label{lem-M=produitfibre}
On a une équivalence
$$\mathcal{S}(\mathrm{M})\simeq\mathcal{S}(\mathbb{R},\mathrm{M})\times_{B_{\mathbb{R}}}\mathcal{T},$$
et un morphisme canonique
$\mathcal{S}(\mathrm{M})\rightarrow\mathcal{S}(\mathbb{R}/W_k)$, où $\mathcal{S}(\mathbb{R}/W_k)$ est le gros topos de l'espace topologique $\mathbb{R}/W_k$.
\end{lem}
\begin{f-proof}
L'objet $\mathfrak{f}^*(E_{\mathbb{R}})$ est donné par la seconde
projection $E_{\mathbb{R}}\times y(\mathbb{R};\mathrm{M})\rightarrow y(\mathbb{R};\mathrm{M})$.
On obtient un isomorphisme (cf \cite{SGA4} IV.5.8.3)
$$\mathcal{S}(\mathbb{R};\mathrm{M})/_{\mathfrak{f}^*(E_{\mathbb{R}})}
=(B_{\mathbb{R}}/_{y(\mathbb{R};\mathrm{M})})/_{(E_{\mathbb{R}}\times
y(\mathbb{R};\mathrm{M}))}
\simeq\mathcal{T}/_{y(\mathrm{M})}=:\mathcal{S}(\mathrm{M}).$$ Les
deux carrés commutatifs
\[ \xymatrix{
\mathcal{S}(\mathrm{M})\simeq
\mathcal{S}(\mathbb{R};\mathrm{M})/_{\mathfrak{f}^*(E_{\mathbb{R}})}
\ar[d]_{} \ar[r]^{}
 & \mathcal{S}(\mathbb{R};\mathrm{M})  \ar[d]_{p}   \\
\mathcal{T}/_{y(\mathbb{R}/W_k)}\simeq
(B_{\mathbb{R}}/_{y(\mathbb{R};\mathbb{R}/W_k)})/_{f^*(E_{\mathbb{R}})}\ar[d]\ar[r]^{}&B_{\mathbb{R}}/_{y(\mathbb{R};\mathbb{R}/W_k)}\ar[d]_f\\
\mathcal{T}\simeq
B_{\mathbb{R}}/_{E_{\mathbb{R}}}\ar[r]^{}&B_{\mathbb{R}} }
\]
sont des pull-backs, où les flèches horizontales sont les morphismes
de localisation. Le résultat suit car $\mathcal{S}(\mathbb{R}/W_k)\simeq\mathcal{T}/_{y(\mathbb{R}/W_k)}$.
\end{f-proof}

\subsection{Le feuilletage} Chaque point
$\overline{u}\in\mathbb{R}/W_k$ est une application continue
$\overline{u}:\{*\}\rightarrow\mathbb{R}/W_k$ qui induit un
morphisme de topos induits $\overline{u}:\mathcal{T}/_{\{*\}}=\mathcal{T}\rightarrow\mathcal{T}/_{y(\mathbb{R}/W_k)}$.
Cette flèche est d'ailleurs une section du morphisme de localisation
$\mathcal{T}/_{y(\mathbb{R}/W_k)}\rightarrow\mathcal{T}$. On définit
la feuille de $\mathcal{S}(\mathrm{M})$ au-dessus de $\overline{u}$
en posant
$$\mathcal{S}(\mathrm{F}_{\overline{u}}):=\mathcal{S}(\mathrm{M})\times_{\mathcal{T}/_{y(\mathbb{R}/W_k)}}\mathcal{T}.$$
Ainsi $\mathcal{S}(\mathrm{F}_{\overline{u}})$ est l'image inverse
par le morphisme
$\mathcal{S}(\mathrm{M})\rightarrow\mathcal{T}/_{\mathbb{R}/W_k}$ du
sous-topos de $\mathcal{T}/_{y(\mathbb{R}/W_k)}$ donné par l'image
du morphisme
$\overline{u}:\mathcal{T}\rightarrow\mathcal{T}/_{y(\mathbb{R}/W_k)}$.
Le lemme \ref{plongementferme-nonequi}, ce morphisme $\overline{u}$ est un plongement fermé. Donc le morphisme
$$\mathcal{S}(\mathrm{F}_{\overline{u}})=\mathcal{S}(\mathrm{M})\times_{\mathcal{T}/_{y(\mathbb{R}/W_k)}}\mathcal{T}
\longrightarrow\mathcal{S}(\mathrm{M})$$ est lui aussi un plongement
fermé. Plus directement, ce morphisme est
induit par l'inclusion fermée de la feuille
$\mathrm{F}_{\overline{u}}\rightarrow\mathrm{M}$. C'est donc un plongement fermé d'après le lemme \ref{plongementferme-nonequi}.

\section{Propriétés analogues du topos
Weil-étale}\label{sect-proprite-analogue}

Soit $Y$ un schéma régulier, séparé et de type fini sur un corps
fini $k=\mathbb{F}_q$. On note $Y_{et}$ le topos étale du schéma $Y$. Le (petit) topos Weil-étale $Y^{sm}_W$ est défini comme la catégorie des faisceaux étales sur $\overline{Y}$ munis d'une action de $W_k$ compatible à celle définie sur $\overline{Y}$ (cf. \cite{Lichtenbaum-finite-field}). D'après \cite{these} Théorème 8.5 (voir aussi \cite{Flach-moi}), on a une équivalence
$$Y^{sm}_W\simeq Y_{et}\times_{B^{sm}_{G_k}}B^{sm}_{W_k}$$
où $G_k$ et $W_k$ désignent le groupe de Galois et le groupe de Weil
de $k$ respectivement. Rappelons que $G_k\simeq\hat{\mathbb{Z}}$ est topologiquement engendré par le Frobenius $\psi$ et que $W_k\simeq\mathbb{Z}$ est le groupe des puissances entières de $\psi$. De plus, le petit topos classifiant $B^{sm}_{G_k}$ (resp. $B^{sm}_{W_k}$) désigne la catégorie des ensembles sur lesquels le groupe $G_k$ (resp. $W_k$) opère continûment. On définit le gros topos classifiant $B_{W_k}$ du groupe discret $W_k$ comme dans \ref{defn-topos-class}.
\begin{defn} Le \emph{gros topos Weil-étale} de $Y$ est défini
par le produit fibré
$$Y_W:= Y_{et}\times_{B^{sm}_{G_k}}B_{W_k}\simeq Y^{sm}_{W}\times_{\underline{Set}}\mathcal{T}.$$
\end{defn}
\subsection{Le morphisme flot et les orbites fermées}
La seconde projection définit un morphisme
$$p:Y_W:=Y_{et}\times_{B^{sm}_{G_k}}B_{W_k}\longrightarrow B_{W_k}.$$
On a aussi un morphisme $B_l:B_{W_k}\rightarrow B_{\mathbb{R}}$
induit par la flèche $l:W_k\rightarrow\mathbb{R}$.
\begin{defn}
Le \emph{morphisme flot} est défini par composition
$$\mathfrak{f}:=B_l\circ p:
Y_W\longrightarrow B_{W_k}\simeq
B_{\mathbb{R}}/_{y(\mathbb{R},\mathbb{R}/W_k)}\longrightarrow
B_{\mathbb{R}}.$$
\end{defn}
Le morphisme $\mathfrak{f}$ se factorise à travers le morphisme de
localisation
$f:B_{\mathbb{R}}/_{y(\mathbb{R},\mathbb{R}/W_k)}\rightarrow
B_{\mathbb{R}}$ et $p$ devient $$p:Y_W\longrightarrow
B_{W_k}\simeq B_{\mathbb{R}}/_{y(\mathbb{R},\mathbb{R}/W_k)}.$$ Un
point fermé $v$ de $Y$ induit un plongement fermé $i_v:B_{W_{k(v)}}\rightarrow Y_W$ (cf. \cite{these} Corollaire 8.20). On a la proposition suivante.
\begin{prop} Un point fermé $v$ de $Y$ induit un plongement fermé
$$i_{v}:B_{W_{k(v)}}\longrightarrow Y_W$$
tel que la composition
$\mathfrak{f}\circ
i_{v}:B_{W_{k(v)}}\rightarrow B_{\mathbb{R}}$ est induite par
le morphisme canonique $W_{k(v)}\rightarrow\mathbb{R}$.
\end{prop}

\subsection{Le topos $\mathcal{Y}$}

On reprend les notations de la section 2.1.2. D'après ce qui précède
$\mathfrak{f}:Y_W\rightarrow B_{\mathbb{R}}$ peut être
vu intuitivement comme le topos
$\mathcal{S}(\mathbb{R};\mathrm{M})\rightarrow B_{\mathbb{R}}$ des
gros $\mathbb{R}$-faisceaux sur $\mathrm{M}$. Pour obtenir
l'analogue du topos $\mathcal{S}(\mathrm{M})$, il suffit de
considérer le pull-back du morphisme $\mathfrak{f}$ le long du
$\mathcal{T}$-point canonique de $B_{\mathbb{R}}$
$$\mathcal{T}\simeq B_{\mathbb{R}}/_{E_{\mathbb{R}}}\rightarrow B_{\mathbb{R}}.$$
\begin{defn}
On définit le topos $\mathcal{Y}$ par le produit fibré
$\mathcal{Y}:=Y_W\times_{B_{\mathbb{R}}}\mathcal{T}$.
\end{defn}
La proposition suivante montre que $\mathcal{Y}$ peut être vu comme un "espace" au-dessus du cercle $\mathbb{R}/W_k$.
\begin{prop}\label{map-EY-dans-S1}
On a un morphisme canonique
$\mathcal{Y}\rightarrow\mathcal{S}(\mathbb{R}/W_k)$.
\end{prop}
\begin{f-proof}
Le diagramme suivant, dans lequel les flèches horizontales sont les
morphismes de localisation, est composé de deux pull-backs.
\[ \xymatrix{
\mathcal{Y}\simeq Y_W/_{\mathfrak{f}^*(E_{\mathbb{R}})}
\ar[d]_{} \ar[r]^{}
 & Y_W  \ar[d]_{p}   \\
\mathcal{T}/_{y(\mathbb{R}/W_k)}\simeq
(B_{\mathbb{R}}/_{y(\mathbb{R};\mathbb{R}/W_k)})/_{f^*(E_{\mathbb{R}})}\ar[d]\ar[r]^{}&B_{\mathbb{R}}/_{y(\mathbb{R};\mathbb{R}/W_k)}\ar[d]_{f}\\
\mathcal{T}\simeq
B_{\mathbb{R}}/_{E_{\mathbb{R}}}\ar[r]^{}&B_{\mathbb{R}} }
\]
La projection $\mathcal{Y}\rightarrow\mathcal{T}$ se factorise donc
à travers le morphisme de localisation
$\mathcal{T}/_{y(\mathbb{R}/W_k)}\rightarrow\mathcal{T}$, pour
induire la flèche
$\mathcal{Y}\rightarrow\mathcal{T}/_{y(\mathbb{R}/W_k)}\simeq\mathcal{S}(\mathbb{R}/W_k)$
\end{f-proof}

\subsection{Le feuilletage}
A nouveau grâce au lemme \ref{plongementferme-nonequi}, chaque point $\overline{u}\in\mathbb{R}/W_k$ définit un
plongement fermé de topos
$$\overline{u}:\mathcal{T}\rightarrow\mathcal{T}/_{y(\mathbb{R}/W_k)}\simeq\mathcal{S}(\mathbb{R}/W_k).$$
\begin{defn}
La feuille $\mathcal{F}_{\overline{u}}$ de $\mathcal{Y}$ au-dessus
de $\overline{u}$ est définie par le produit fibré
$$\mathcal{F}_{\overline{u}}:=\mathcal{Y}\times_{\mathcal{S}(\mathbb{R}/W_k)}\mathcal{T}.$$
\end{defn}
La seconde projection définit donc un morphisme canonique
$$\mathcal{F}_{\overline{u}}=\mathcal{Y}\times_{\mathcal{S}(\mathbb{R}/W_k)}\mathcal{T}
\longrightarrow\mathcal{Y}.$$
\begin{prop} Le morphisme précédent
$\mathcal{F}_{\overline{u}}\rightarrow\mathcal{Y}$ est un plongement
fermé.
\end{prop}
\begin{f-proof}
Par définition, $\mathcal{F}_{\overline{u}}$ est l'image inverse par
le morphisme
$\mathcal{Y}\rightarrow\mathcal{S}(\mathbb{R}/W_k)$ du
sous-topos fermé
$\overline{u}:\mathcal{T}\rightarrow\mathcal{S}(\mathbb{R}/W_k)$.
Il suit que
$$\mathcal{F}_{\overline{u}}=\mathcal{Y}\times_{\mathcal{S}(\mathbb{R}/W_k)}\mathcal{T}
\longrightarrow\mathcal{Y}$$ est un plongement fermé.
\end{f-proof}
\begin{prop}\label{lemfeuille} Soit $\overline{u}$ un point de $\mathbb{R}/W_k$. On a
une équivalence canonique
$$\mathcal{F}_{\overline{u}}\simeq
\overline{Y}_{et}\times\mathcal{T},$$ où
$\overline{Y}_{et}$ est le topos étale du schéma
$\overline{Y}=Y\times_k\overline{k}$.
\end{prop}
\begin{f-proof}
La preuve sera donnée dans la démonstration du théorème \ref{propmorphisme} (iii).
\end{f-proof}

Les topos $\overline{Y}_{et}\times\mathcal{T}$ et $\overline{Y}_{et}$ sont cohomologiquement équivalents (cf. \cite{these}). Il est par ailleurs bien connu que $\overline{Y}_{et}$ est de dimension cohomologique $2.dim(Y)$.

\begin{rem}Soit $Y$ un schéma régulier de dimension $d$ sur
un corps fini $k$. Alors le gros topos Weil-étale $Y_W$ peut
être vu comme l'analogue du topos
$\mathcal{S}(\mathbb{R};\mathrm{M})$, où $\mathrm{M}$ est un espace de dimension $2d+1$. Les feuilles
$$\mathcal{F}_{\overline{u}}\simeq
\overline{Y}_{et}\times\mathcal{T}\longrightarrow\mathcal{Y}$$
de $\mathcal{Y}$ sont analogues à des sous-espaces fermés de
codimension $1$ dans
$\mathrm{M}$.
\end{rem}

\subsection{Le morphisme du topos dynamique dans le topos
Weil-étale}\label{sectionmorphisme}

Soit $Y$ un schéma lisse et propre sur un corps fini $k=\mathbb{F}_q$. Supposons que l'on puisse associer fonctoriellement à
$Y$ un système dynamique feuilleté
$$(\mathrm{M}_{Y},\mathrm{F},\phi)$$
de la forme décrite dans l'introduction de la section \ref{sect-systdyn-charp} (cf.
\cite{Deninger-explicit formulas}, \cite{Deninger-NTandDSonFS},
\cite{Deninger-possible-significance}, \cite{Leichtnam-invitation} Section 4.2 et \cite{Leichtnam} Open Question 2). On note
$\mathfrak{D}(Y)=(\mathbb{R},\mathrm{M}_{Y})$ l'action continue de
$\mathbb{R}$ sur l'espace topologique $\mathrm{M}_{Y}$. On note
$\gamma_v:\mathbb{M}_{N(v)}\rightarrow\mathrm{M}_{Y}$ l'orbite fermée dans $\mathrm{M}_{Y}$ de longueur
$log(N(v))$ correspondant, via le foncteur $\mathfrak{D}$, à l'inclusion d'un point fermé $Spec(k(v))\rightarrow Y$.
Il existe alors un foncteur
$$\mathfrak{D}:Et_Y\longrightarrow {Top}^{\mathbb{R}}/_{(\mathbb{R},\mathrm{M}_{Y})},$$
où $Et_Y$ désigne la catégorie des schémas étales au-dessus de $Y$.
\begin{hypothese}\label{hypo-charp}
Nous supposons que les propriétés suivantes sont satisfaites.
\begin{itemize}
\item Le foncteur $\mathfrak{D}$ est exact à gauche (i.e. il commute aux produits fibrés).
\item Un morphisme étale $U\rightarrow V$ au-dessus de $Y$ est envoyé sur
(une application continue $\mathbb{R}$-équivariante qui est) un
étalement.
\item Une famille surjective de morphismes étales de schémas est envoyée sur une famille surjective d'applications continues.
\end{itemize}
\end{hypothese}

\begin{thm}\label{propmorphisme}
Si le foncteur $\mathfrak{D}$ existe et satisfait l'hypothèse \ref{hypo-charp}, alors on a les résultats suivants.

\emph{\,\,(i)} Il existe un morphisme canonique
$$d:\mathcal{S}(\mathbb{R};\mathrm{M}_{Y})\longrightarrow Y_W$$ défini au-dessus de $B_{W_{k}}$. Ce morphisme induit une flèche $$\mathcal{S}(\mathrm{M}_{Y})\longrightarrow {Y_W}\times_{B_{\mathbb{R}}}\mathcal{T}=:\mathcal{Y}$$
définie au-dessus de $\mathcal{S}(\mathbb{R}/W_k)$.

\emph{\,\,(ii)} Soit $v$ un point fermé de $Y$. Le diagramme suivant est un pull-back.
\[ \xymatrix{
 B_{W_{k(v)}}\ar[r]^{Id}\ar[d]^{i_{\gamma_v}}& B_{W_{k(v)}}\ar[d]^{i_v}
 \\
  \mathcal{S}(\mathbb{R};\mathrm{M}_{Y})\ar[r]^{d}&Y_W
} \]

\emph{\,\,(iii)} Soient $\overline{u}$ un point du cercle
$\mathbb{R}/W_{\mathbb{F}_q}$. On note
$\mathcal{S}(\mathrm{F}_{\overline{u}})\rightarrow
\mathcal{S}(\mathrm{M}_Y)$ et
$\mathcal{F}_{\overline{u}}\rightarrow\mathcal{Y}$ les fibres des topos $\mathcal{S}(\mathrm{M}_Y)$ et $\mathcal{Y}$ au-dessus de $\overline{u}\in \mathbb{R}/W_{\mathbb{F}_q}$. Alors les feuilles
$\mathcal{S}(\mathrm{F}_{\overline{u}})$ et $\mathcal{F}_{\overline{u}}$ des topos $\mathcal{S}(\mathrm{M}_{Y})$ et $\mathcal{Y}$ sont des sous-topos fermés, et $\mathcal{S}(\mathrm{F}_{\overline{u}})$ est l'image inverse de $\mathcal{F}_{\overline{u}}$ à travers le morphisme
$\mathcal{S}(\mathrm{M}_{Y})
\rightarrow\mathcal{Y}$. Enfin, on a une équivalence
$$\mathcal{F}_{\overline{u}}\simeq\overline{Y}_{et}\times\mathcal{T}.$$
\end{thm}
\begin{f-proof}

(i) Puisqu'un étalement admet des sections locales au-dessus de son
image, on a un morphisme de sites exacts à gauche
$$\mathfrak{D}:(Et_Y;\mathcal{J}_{et})\longrightarrow({Top}^{\mathbb{R}}/_{(\mathbb{R},\mathrm{M}_{Y})},\mathcal{J}_{ls}),$$
et donc un morphisme de topos
$\mathcal{S}(\mathbb{R};\mathrm{M}_{Y})\rightarrow Y_{et}$.
On considère maintenant le diagramme de sites exacts à gauche suivant.
\[ \xymatrix{
 (Et_Y;\mathcal{J}_{et})  \ar[r]
 &({Top}^{\mathbb{R}}/_{(\mathbb{R},\mathrm{M}_{Y})};\mathcal{J}_{ls})    \\
 (Et_{Spec(\mathbb{F}_q)};\mathcal{J}_{et})\ar[r]\ar[u]
 &({Top}^{\mathbb{R}}/_{\mathbb{M}_{q}};\mathcal{J}_{ls})\ar[u]
} \] Ci-dessus, les flèches verticales sont données par changements
de bases. Ce diagramme est (pseudo)-commutatif, car le foncteur $\mathfrak{D}$ a été supposé exact à
gauche. On en déduit un diagramme commutatif de topos
\[ \xymatrix{
 \mathcal{S}(\mathbb{R};\mathrm{M}_{Y})\ar[r]^{}\ar[d]&Y_{et}\ar[d]   \\
 B_{W_{k}}\simeq B_{\mathbb{R}}/_{y(\mathbb{M}_q)}\ar[r]& B^{sm}_{G_{k}}
} \] Par la propriété universelle du
produit fibré, on obtient un morphisme
\begin{equation}\label{morphisme-d}
d:\mathcal{S}(\mathbb{R};\mathrm{M}_{Y})\longrightarrow Y_W:=Y_{et}\times_{B^{sm}_{G_{k}}}B_{W_{k}}
\end{equation}
tel que le diagramme suivant soit commutatif.
\[ \xymatrix{
 \mathcal{S}(\mathbb{R};\mathrm{M}_{Y})\ar[r]^{d}\ar[d]&Y_W\ar[d]   \\
 B_{\mathbb{R}}/_{y(\mathbb{M}_q)}\ar[r]^{\simeq}  &B_{W_{k}}
} \]
Le morphisme $d$ ci-dessus est donc défini au-dessus de $B_{\mathbb{R}}$. Soit $\mathcal{T}\rightarrow B_{\mathbb{R}}$ le $\mathcal{T}$-point
canonique de $B_{\mathbb{R}}$. D'après le lemme \ref{lem-M=produitfibre}, on obtient un morphisme
$$(d,Id_{\mathcal{T}}):\mathcal{S}(\mathrm{M}_{Y})\simeq
\mathcal{S}(\mathbb{R};\mathrm{M}_{Y})\times_{B_{\mathbb{R}}}\mathcal{T}
\longrightarrow {Y_W}\times_{B_{\mathbb{R}}}\mathcal{T}=:\mathcal{Y}.$$
défini au-dessus de $B_{W_k}\times_{B_{\mathbb{R}}}\mathcal{T}\simeq\mathcal{S}(\mathbb{R}/W_k)$.

(ii) Soit $v$ un point fermé de $Y$. Il lui correspond une orbite fermée
$\gamma_v$ dans $\mathrm{M}_Y$ de longueur $log(N(v))$, où $N(v)$
est le cardinal du corps résiduel $k(v)$. Le morphisme de topos
$$i_{\gamma_v}:B_{W_{k(v)}}\longrightarrow\mathcal{S}(\mathbb{R};\mathrm{M}_{Y})$$
défini dans la proposition \ref{orbite-ferme-syst-dyn-char-p} est
induit par le morphisme de sites exacts à gauche
$$
\fonc{i_{\gamma_v}^*}{(Top^{\mathbb{R}}/_{(\mathbb{R},\mathrm{M}_{Y})};\mathcal{J}_{ls})}
{(Top^{\mathbb{R}}/_{\mathbb{M}_{N(v)}};\mathcal{J}_{ls})}
{Z}{Z\times_{\mathrm{M}_{Y}}\mathbb{M}_{N(v)}}
$$
où $\mathbb{R}$ opère diagonalement sur
$Z\times_{\mathrm{M}_{Y}}\mathbb{M}_{N(v)}$. Le diagramme
suivant de sites exacts à gauche est commutatif car $\mathfrak{D}$ commute aux produits fibrés.
\[ \xymatrix{
 {Top}^{\mathbb{R}}/_{\mathbb{M}_{N(v)}}  &\ar[l]_{\alpha_v^*} Et_{Spec(k(v))}   \\
 {Top}^{\mathbb{R}}/_{(\mathbb{R},\mathrm{M}_{Y})}\ar[u]_{i_{\gamma_v}^*}
 & \ar[l]_{\mathfrak{D}}\ar[u]_{u_v^*}Et_{Y}
} \] Ce diagramme de sites induit donc un diagramme commutatif de topos. De plus,
l'image de $B^{sm}_{G_{k(v)}}$ dans $Y_{et}$ est le
sous-topos fermé complémentaire de l'ouvert défini par l'objet
$y(U)=y(Y-\{v\})$ de $Y_{et}$, qui est un sous-objet de
l'objet final. De la même manière, le lemme \ref{plongementferme-equi} montre que l'image de $B_{W_{k(v)}}$ dans
$\mathcal{S}(\mathbb{R};\mathrm{M}_{Y})$ est le sous-topos fermé
complémentaire de l'ouvert défini par l'objet
$\gamma^*(y(U))=y(\mathrm{M}_{U})$ de
$\mathcal{S}(\mathbb{R};\mathrm{M}_{Y})$. Il suit que le diagramme suivant est un pull-back.
\[ \xymatrix{
 B_{W_{k(v)}} \ar[d]_{i_{\gamma_v}} \ar[r]^{\alpha_v} &B^{sm}_{G_{k(v)}} \ar[d]_{u_v}   \\
 \mathcal{S}(\mathbb{R};\mathrm{M}_{Y})\ar[r]&Y_{et}
} \]
Considérons maintenant le diagramme commutatif
suivant.
\[ \xymatrix{
 B_{W_{k(v)}}\ar[r]^{Id}\ar[d]^{i_{\gamma_v}}& B_{W_{k(v)}}\ar[d]^{i_v}\ar[r]^{\alpha_v}& B^{sm}_{G_{k(v)}}\ar[d]^{u_v}\\
  \mathcal{S}(\mathbb{R};\mathrm{M}_{Y})\ar[r]^{d}&Y_W
  \ar[r]
  &{Y}_{et}
} \] Ici, $\gamma_v$ désigne l'orbite fermée de $\mathrm{M}_{Y}$
correspondant au point fermé $v\in Y$.  Le carré de
droite est un pull-back et nous venons de montrer que le carré total est un pull-back. Il suit que le carré de
gauche est lui aussi un pull-back.

(iii) Le morphisme $(d,Id_{\mathcal{T}})$ est défini au-dessus de
$B_{\mathbb{R}}/_{y(\mathbb{M}_q)}\times_{B_{\mathbb{R}}}\mathcal{T}\simeq\mathcal{S}(\mathbb{R}/W_{\mathbb{F}_q})$. Soient $\overline{u}$ un point du cercle
$\mathbb{R}/W_{\mathbb{F}_q}$,
$\overline{u}:\mathcal{T}\rightarrow\mathcal{S}(\mathbb{R}/W_{\mathbb{F}_q})$
le $\mathcal{T}$-point de $\mathcal{S}(\mathbb{R}/W_{\mathbb{F}_q})$ correspondant et
$\mathrm{F}_{\overline{u}}$ la feuille de $\mathrm{M}_{Y}$ au-dessus
de $\overline{u}$. Alors on a
$$\mathcal{S}(\mathrm{F}_{\overline{u}})
=\mathcal{S}(\mathrm{M}_Y)\times_{\mathcal{S}(\mathbb{R}/W_{\mathbb{F}_q})}\mathcal{T}
=\mathcal{S}(\mathrm{M}_Y)\times_{\mathcal{Y}}\mathcal{Y}
\times_{\mathcal{S}(\mathbb{R}/W_{\mathbb{F}_q})}\mathcal{T}
=\mathcal{S}(\mathrm{M}_Y)\times_{\mathcal{Y}}\mathcal{F}_{\overline{u}},$$
où $\mathcal{F}_{\overline{u}}$ est la feuille de $\mathcal{Y}$
au-dessus de $\overline{u}$.  En d'autres termes, le diagramme suivant est composé de pull-backs :
\[ \xymatrix{
\mathcal{S}(\mathrm{F}_{\overline{u}})\ar[r]\ar[d]&\mathcal{F}_{\overline{u}}\ar[d]\ar[r]&\mathcal{T}\ar[d]^{\overline{u}}   \\
 \mathcal{S}(\mathrm{M}_{Y})\ar[r]^{\,\,\,\,\,\,(d,Id_{\mathcal{T}})}  &\mathcal{Y}\ar[r]&\mathcal{S}(\mathbb{R}/W_{\mathbb{F}_q})
} \]
D'après le lemme \ref{plongementferme-nonequi}, le morphisme $\overline{u}:\mathcal{T}\rightarrow\mathcal{S}(\mathbb{R}/W_{\mathbb{F}_q})$ est un plongement fermé. Il suit que les morphismes
$\mathcal{S}(\mathrm{F}_{\overline{u}})\rightarrow
\mathcal{S}(\mathrm{M}_Y)$ et
$\mathcal{F}_{\overline{u}}\rightarrow\mathcal{Y}$ sont des plongements fermés.

Il reste à définir l'équivalence $\mathcal{F}_{\overline{u}}\simeq\overline{Y}_{et}\times\mathcal{T}$. Les deux carrés commutatifs suivants sont des pull-backs.
\[ \xymatrix{
\mathcal{F}_{\overline{u}}\ar[r]\ar[d] & \mathcal{Y}\simeq
Y_W/_{\mathfrak{f}^*(E_{\mathbb{R}})} \ar[d]_{}
\ar[r]^{}
 & Y_W  \ar[d]_{p}   \\
\mathcal{T}\ar[r]&\mathcal{S}(\mathbb{R}/W_k)\simeq
(B_{\mathbb{R}}/_{y(\mathbb{R};\mathbb{R}/W_k)})/_{f^*(E_{\mathbb{R}})}\ar[r]^{}
&B_{\mathbb{R}}/_{y(\mathbb{R};\mathbb{R}/W_k)}}
\]
Le carré total est aussi un pull-back. Il s'identifie à
\[ \xymatrix{
\mathcal{F}_{u}\ar[r]\ar[d]
 & Y_{et}\times_{B^{sm}_{G_k}}B_{W_k}  \ar[d]_{p}   \\
\mathcal{T}\ar[r] &B_{W_k} }\] où la flèche $\mathcal{T}\rightarrow
B_{W_k}$ est le $\mathcal{T}$-point canonique de $B_{W_k}$. On a
donc
\begin{equation}\label{feuille-identific-0}
\mathcal{F}_{u}\simeq Y_{et}\times_{B^{sm}_{G_k}}B_{W_k}\times_{B_{W_k}}\mathcal{T}
\simeq Y_{et}\times_{B^{sm}_{G_k}}\mathcal{T}.
\end{equation}
Mais le morphisme $\mathcal{T}\rightarrow B^{sm}_{G_k}$ se factorise
à travers $\underline{Set}\rightarrow B^{sm}_{G_k}$, le point
canonique de $B^{sm}_{G_k}$. On en déduit l'identification
\begin{equation}\label{feuille-identific-1}
Y_{et}\times_{B^{sm}_{G_k}}\mathcal{T}
\simeq Y_{et}\times_{B^{sm}_{G_k}}\underline{Set}\times_{\underline{Set}}\mathcal{T}.
\end{equation}
Montrons l'équivalence suivante :
\begin{equation}\label{feulle-identific-2}
\overline{Y}_{et}
\simeq Y_{et}\times_{B^{sm}_{G_k}}\underline{Set}.
\end{equation}
On a les équivalences
$$\underline{Set}=Spec(\overline{k})_{et}=\underleftarrow{lim}\,Spec(k^n)_{et}=\underleftarrow{lim}\,B^{sm}_{G_{k^n}},$$ où $k^n/k$ est l'unique extension de degré $n$ (dans une clôture algébrique fixée). Par ailleurs, le morphisme
$B^{sm}_{G_{k^n}}\rightarrow B^{sm}_{G_{k}}$ induit par l'inclusion
$G_{k^n}\rightarrow G_{k}$, s'identifie au morphisme de localisation
$$B^{sm}_{G_{k}}/_{(G_{k}/G_{k^n})}\longrightarrow B^{sm}_{G_{k}}.$$
Le pull-back de l'objet $G_{k}/G_{k^n}$ de $B^{sm}_{G_{k}}$ par le
morphisme $Y_{et}\rightarrow{B^{sm}_{G_k}}$ est
précisément l'objet de $Y_{et}$ représenté par
$Y^n:=Y\times_kk^n$. On en déduit
$$Y_{et}\times_{B^{sm}_{G_k}}B^{sm}_{G_{k^n}}=Y_{et}\times_{B^{sm}_{G_k}}B^{sm}_{G_{k}}/_{(G_{k}/G_{k^n})}
=Y_{et}/_{y(Y^n)}=Y^n_{et}.$$ On obtient
(\ref{feulle-identific-2}) grâce aux identifications
$$Y_{et}\times_{B^{sm}_{G_k}}\underline{Set}=
\underleftarrow{lim}\,Y_{et}\times_{B^{sm}_{G_k}}B^{sm}_{G_{k^n}}
=\underleftarrow{lim}\,Y^n_{et}=\overline{Y}_{et}.$$
Ci-dessus, on utilise le fait que la limite projective des topos
étales d'une famille filtrante de schémas $Y^n$ (respectivement
$Spec(k^n)$) quasi-séparés, quasi-compacts et dont les morphismes de
transition sont affines, s'identifie au topos étale du schéma limite
projective $\overline{Y}$ (respectivement $Spec(\overline{k})$). Les
équivalences (\ref{feuille-identific-0}), (\ref{feuille-identific-1}) et
(\ref{feulle-identific-2}) donnent
$$\mathcal{F}_{\overline{u}}\simeq\overline{Y}_{et}\times\mathcal{T}.$$

\end{f-proof}

L'existence de $\mathrm{M}_{Y}$ est suggérée dans \cite{Leichtnam-invitation} Section 4.2 pour toute variété $Y$ lisse et projective. Cependant, on s'attend à ce que la flèche $Y\rightarrow\mathrm{M}_{Y} $ ne soit fonctorielle qu'en se restreignant aux variétés ordinaires (cf. \cite{Deninger-possible-significance} Section 4.7 et \cite{Leichtnam-invitation} Section 4.2 Remark 1) (voir aussi \cite{Milne-motivesfinitefields} Theorem 2.49). Le lecteur soucieux de ce problème pourra donc supposer que $Y$ est ordinaire, restreindre $Et_Y$ à ces variétés et munir la sous-catégorie pleine de $Et_Y$ obtenue de la topologie induite. Le résultat précédent est alors valable à cette modification du topos étale $Y_{et}$ près.

Nous ignorons ce problème dans cette section car les résultats présentés dans cet article suggèrent que le système dynamique de Deninger n'existe pas sous la forme d'un espace au sens classique (par exemple à cause de l'argument de  \cite{Deninger-possible-significance} Section 4.7, ou encore comme le suggère la structure topologique du groupe fondamental \cite{Fund-group-I}) mais sous la forme plus générale d'un topos. Ce topos étant par définition la structure topologique sous-jacente à la cohomologie de \cite{Deninger4} Conjecture 2, il doit exister (fonctoriellement) pour tout schéma séparé de type fini au-dessus de $Spec(\mathbb{Z})$. Le théorème précédent établit le lien existant entre ce topos conjectural et le topos Weil-étale.

\section{Le système dynamique conjecturalement associé à l'anneau d'entiers d'un corps de nombres}

C. Deninger suggère l'existence d'un foncteur de la catégorie des
schémas plats séparés de type fini sur $Spec(\mathbb{Z})$ dans celle des
systèmes dynamiques munis d'un feuilletage (cf.
\cite{Deninger-some-analogies},
\cite{Deninger-possible-significance}, \cite{Deninger-NTandDSonFS},
\cite{Deninger2}, \cite{Deninger-explicit formulas}). Précisons que Deninger considère cette conjecture comme optimiste, et qu'il envisage que ce système dynamique puisse prendre la forme d'un site (i.e. d'un topos). La forme conjecturale précise du système dynamique associé à $Spec(\mathbb{Z})$ est donnée dans \cite{Leichtnam-invitation} Conjecture 1.

Un tel morphisme $X\rightarrow Spec(\mathbb{Z})$ induirait donc un
morphisme de systèmes dynamiques. On se restreint dans cette partie aux
anneaux d'entiers de corps de nombres. Soit $K$ un corps de nombres,
$X:=Spec(\mathcal{O}_K)$ le spectre de l'anneau d'entiers de $K$,
$X_{\infty}$ l'ensemble des places archimédiennes de $K$ et
$\bar{X}:=(X;X_{\infty})$ la compactification d'Arakelov de
$X$. Le système dynamique $(\mathrm{M}_{X},\phi)$ conjecturalement associé à $X$ est
muni d'un feuilletage $\mathrm{F}$ de codimension $1$ dont les
feuilles sont partout perpendiculaires aux trajectoires du flot. Alors $\mathrm{M}_{\bar{X}}$
est une compactification de $\mathrm{M}_{X}$, de dimension topologique $3$. L'espace $\mathrm{M}_{X}$ a une structure d'espace laminé (cf. \cite{Deninger-NTandDSonFS} Section 5.1). En particulier, $\mathrm{M}_{X}$ est localement homéomorphe à un produit $U\times T$, où $U$ est un ouvert non-vide de $\mathbb{R}^3$ et $T$ un espace totalement discontinu. Le feuilletage $\mathrm{F}$ est une partition de $\mathrm{M}_{X}$ par des espaces laminés de dimension 2, qui est localement triviale (cf. \cite{Deninger-NTandDSonFS} Section 5.2). Le flot est donné par un morphisme
$$
\fonc{\phi}{\mathbb{R}\times \mathrm{M}_{\bar{X}}}
{\mathrm{M}_{\bar{X}}}{(t;m)}{\phi^t(m)}
$$
respectant les conditions d'une action. Ce flot est compatible au
feuilletage $\mathrm{F}$. Plus précisément, $\phi^t$ induit un
isomorphisme $\mathrm{F}_m\rightarrow\mathrm{F}_{\phi^t(m)}$ quel que soit
le couple $(t;m)$, où $\mathrm{F}_m$ désigne la feuille contenant
$m\in\mathrm{M}$. Un point fermé $v\rightarrow X$ correspond à une
orbite fermée
$$\gamma_v:\mathbb{R}/log(N(v))\mathbb{Z}\longrightarrow\mathrm{M}_{\bar{X}},$$
où $N(v):=\mid k(v)\mid $ est la norme de $v$. Le groupe
$\mathbb{R}$ opère naturellement sur l'espace homogène
$\mathbb{R}/log(N(v))\mathbb{Z}$ et l'immersion $\gamma_v$ est
$\mathbb{R}$-équivariante. On note $l(\gamma_v)=log(N(v))$ la longueur de l'orbite $\gamma_v$. En oubliant le flow et le feuilletage, $\mathrm{M}_{\bar{X}}$ devient un espace de dimension trois dans lequel une place finie de $K$ correspond à un noeud (cf. \cite{moi} et \cite{Morishita-primesandknots}).

Une place archimédienne $\mathfrak{p}\in X_{\infty}\subset\bar{X}$ correspond à
un point fixe du flot $m_{\mathfrak{p}}\in \mathrm{M}_{\bar{X}}$. La feuille $\mathrm{F}_{m_{\mathfrak{p}}}$ passant
par un tel point fixe est globalement stabilisée par le flot.  A nouveau, le topos que nous associons à $\mathrm{M}_{X}$ ne retient que la structure d'espace topologique muni d'une action continue de $\mathbb{R}$.

\subsection{Le flot, les orbites fermées et les points fixes}\label{subsect-orbites-fermee-char-0}
Le topos $\mathcal{S}(\mathbb{R};\mathrm{M}_{\bar{X}})$ est
défini comme le topos induit
$$\mathcal{S}(\mathbb{R},\mathrm{M}_{\bar{X}}):=B_{\mathbb{R}}/_{y(\mathbb{R},\mathrm{M}_{\bar{X}})}.$$
\begin{defn}\label{def-flot}
Le \emph{morphisme flot} est le morphisme de localisation
$$\mathfrak{f}:\mathcal{S}(\mathbb{R};\mathrm{M}_{\bar{X}})\longrightarrow B_{\mathbb{R}}.$$
\end{defn}
Ce morphisme est d'ailleurs induit par l'application continue
$\mathbb{R}$-équivariante $\mathrm{M}_{\bar{X}}\rightarrow\{*\}$, où $\{*\}$ est l'espace ponctuel sur lequel $\mathbb{R}$ opère
trivialement. Le morphisme flot $\mathfrak{f}$ \emph{détermine}
l'action de $\mathbb{R}$ sur $\mathrm{M}_{\bar{X}}$ (cf. Section \ref{appendix}), d'où la terminologie. On considère un point fermé $v$ de $X\subset\bar{X}$. Il lui correspond
une orbite fermée $\gamma_v:\mathbb{R}/log(N(v))\mathbb{Z}\rightarrow\mathrm{M}_{\bar{X}}$.
On considère de plus une place archimédienne $\mathfrak{p}\in X_{\infty}$  du corps
de nombres $K$. On note $m_{\mathfrak{p}}\in\mathrm{M}_{\bar{X}}$ le
point fixe du flot correspondant.

\begin{thm}\label{orbiteferme+ptsfixes}
Il existe un morphisme canonique $\mathfrak{f}:\mathcal{S}(\mathbb{R};\mathrm{M}_{\bar{X}})\rightarrow B_{\mathbb{R}}$.

\emph{(i)} L'orbite fermée $\gamma_v$ induit un plongement fermé
$i_v:B_{W_{k(v)}}\rightarrow\mathcal{S}(\mathbb{R};\mathrm{M}_{\bar{X}})$
tel que la composition $\mathfrak{f}\circ
i_v:B_{W_{k(v)}}\rightarrow B_{\mathbb{R}}$ est induite par le morphisme canonique $W_{k(v)}\rightarrow\mathbb{R}$.

\emph{(ii)} Le point fixe $m_{\mathfrak{p}}$ induit un plongement fermé
$i_{\mathfrak{p}}:B_{\mathbb{R}}\rightarrow\mathcal{S}(\mathbb{R};\mathrm{M}_{\bar{X}})$
tel que la composition
$\mathfrak{f}\circ i_{\mathfrak{p}}:B_{\mathbb{R}}\rightarrow B_{\mathbb{R}}$
est l'identité de $B_{\mathbb{R}}$.
\end{thm}

\begin{f-proof}
Le morphisme $\mathfrak{f}$ est donné par la définition \ref{def-flot}. La flèche $\gamma_v$ induit un morphisme
$\gamma_v: \mathbb{M}_{N(v)}\rightarrow(\mathbb{R};\mathrm{M}_{\bar{X}})$
dans la catégorie ${Top}^{\mathbb{R}}$. En appliquant le foncteur de Yoneda, on obtient un morphisme
$$y(\mathbb{M}_{N(v)})=y(\mathbb{R},\mathbb{R}/log\,N(v))\longrightarrow y(\mathbb{R};\mathrm{M}_{\bar{X}})$$
dans la catégorie
$B_{\mathbb{R}}\simeq\widetilde{({Top}^{\mathbb{R}};\mathcal{J}_{ls})}$
qui induit à son tour un morphisme de topos induits
$$i_{\gamma_v}:B_{W_{k(v)}}\simeq B_{\mathbb{R}}/_{y(\mathbb{M}_{N(v)})}\longrightarrow B_{\mathbb{R}}/_{y(\mathbb{R};\mathrm{M}_{\bar{X}})}.$$
La composition $\mathfrak{f}\circ i_{\gamma_v}$ est donc donnée par le morphisme de
localisation $B_{\mathbb{R}}/_{y(\mathbb{M}_{N(v)})}\rightarrow
B_{\mathbb{R}}$, qui s'identifie à la flèche
$B_{W_{k(v)}}\rightarrow B_{\mathbb{R}}$ induite par le morphisme de
groupes topologiques $l_v:W_{k(v)}\rightarrow\mathbb{R}$ envoyant le
générateur canonique de $W_{k(v)}$ sur $log(N(v))$. Le point fixe $m_{\mathfrak{p}}$ induit une flèche
$m_{\mathfrak{p}}: \{*\}\rightarrow(\mathbb{R},\mathrm{M}_{\bar{X}})$
dans la catégorie ${Top}^{\mathbb{R}}$. On obtient un morphisme de topos
$$i_{\mathfrak{p}}:B_{\mathbb{R}}=\mathcal{S}(\mathbb{R};*)\longrightarrow\mathcal{S}(\mathbb{R};\mathrm{M}_{\bar{X}}).$$
Le morphisme composé $\mathfrak{f}\circ i_{\mathfrak{p}}$ est induit par la
composition
$\{*\}\rightarrow(\mathbb{R},\mathrm{M}_{\bar{X}})\rightarrow\{*\}$,
qui est l'identité du point dans la catégorie
${Top}^{\mathbb{R}}$. Ainsi, on a un isomorphisme
canonique $$\mathfrak{f}\circ i_{\mathfrak{p}}\simeq Id_{B_{\mathbb{R}}}.$$
Il reste à montrer que les morphismes $i_{\gamma_v}$ et $i_{\mathfrak{p}}$ sont des plongements fermés. Les morphismes $i_{\gamma_v}$ et $i_{\mathfrak{p}}$ sont des morphismes de localisation associés à des sous-espaces fermés stables sous l'action de $\mathbb{R}$. Les deux lemmes suivants montrent qu'un tel morphisme est un plongement fermé.
\begin{lem}\label{plongementferme-nonequi}
Soit $i:Z\rightarrow M$ un sous-espace espace fermé d'un espace topologique $M$.
Alors le morphisme de topos $\mathcal{S}(Z)\rightarrow\mathcal{S}(M)$ est le plongement fermé complémentaire du sous-topos ouvert $\mathcal{S}(U)\rightarrow\mathcal{S}(M)$.
\end{lem}
On considère la catégorie $TOP_M$ des systèmes $$(F_{M'},\varphi_f)_{M'\in Ob(Top/M),\,f\in Fl(Top/M)}$$
définis comme suit. Pour tout espace $M'$ au-dessus de $M$, $F_{M'}$ est un petit faisceau (i.e. un espace étalé) sur $M'$. Pour une fonction continue $f:M''\rightarrow M'$ au-dessus de $M$, $\varphi_f:f^*F_{M'}\rightarrow F_{M''}$ est un morphisme de faisceaux sur $M''$. Les morphismes $\varphi_f$ satisfont la condition de transitivité $\varphi_{f\circ g}=\varphi_g\circ g^*\varphi_f$, et $\varphi_f$ est un isomorphisme lorsque $f$ est un étalement. D'après (\cite{SGA4} IV.4.10.3), la catégorie $TOP_M$ est naturellement équivalente au gros topos $\mathcal{S}(M)\simeq TOP(M)$.

Les inclusions du sous-espace fermé $Z\hookrightarrow M$ et de son complémentaire ouvert $U\hookrightarrow M$ induisent des morphismes de topos
$$\underline{i}:TOP_Z\longrightarrow TOP_M \longleftarrow TOP_U:\underline{j}$$
Le morphisme $\underline{j}$ est le plongement ouvert donné par le sous-objet $(M'_U,\varphi_f)$ de l'objet final, où le faisceau $M'_U$ sur $M'$ est défini par l'espace étalé $M'\times_M U\rightarrow M'$ et $\varphi_f=Id_{M''\times_MU}$. Le foncteur $\underline{j}^*$ peut être décrit comme suit : $\underline{j}^*(F_{M'},\varphi_f)=(F_{U'},\varphi_g)$, où $U'\rightarrow U$ est vu comme un espace au-dessus de $M$ et $g:U''\rightarrow U'$ comme une flèche au-dessus de $M$.  Le foncteur $\underline{i}^*$ peut être décrit de manière analogue, i.e. en écrivant $\underline{i}^*(F_{M'},\varphi_f)=(F_{Z'},\varphi_g)$. Il est clair que les foncteurs $\underline{j}^*$ et $\underline{i}^*$ commutent aux limites projectives et inductives quelquonques, car ces limites se calculent argument par argument. On a donc bien deux morphismes (essentiels) de topos $\underline{i}$ et $\underline{j}$ (en fait $\underline{i}$ et $\underline{j}$ sont des plongements ouverts).

Le foncteur $\underline{i}_*$ est défini de la manière suivante. Soit $(L_{Z'},\phi_g)$ un objet de $TOP_Z$. Pour $M'\rightarrow M$, on note $i':Z'=M'\times_MZ\rightarrow M'$ et on définit $F_{M'}:=i'_*L_{Z'}$. Si $f:M''\rightarrow M'$ est une fonction continue au-dessus de $M$, alors $g:Z''\rightarrow Z'$ est définie par changement de base. On a donc un morphisme $\phi_g:g^*L_{Z'}\rightarrow L_{Z''}$. En appliquant le foncteur $i''_*$, on obtient un morphisme $i''_*(\phi_g):i''_*g^*L_{Z'}\rightarrow i''_*L_{Z''}$. Considérons le carré cartésien :
\[ \xymatrix{
Z'' \ar[d]_g \ar[r]^{i''} &M'' \ar[d]_f   \\
Z'\ar[r]^{i'}&M'
} \]
Comme $i'$ est une immersion fermée, on a $f^*i'_*\simeq i''_*g^*$, et on obtient le morphisme
$$\varphi_f:f^*i'_*L_{Z'}\simeq i''_*g^*L_{Z'}\rightarrow i''_*L_{Z''}$$
qui est d'ailleurs un isomorphisme lorsque $f$ est un étalement. On définit ainsi $\underline{i}_*(L_{Z'},\phi_g)=(F_{M'},\varphi_f)$ avec $F_{M'}:=i'_*L_{Z'}$ et $\varphi_f$ défini ci-dessus. On montre facilement la formule d'adjonction
$$Hom_{TOP_Z}(\underline{i}^*(\mathcal{F}_{M'},\varphi_f),(L_{Z'},\phi_g))=Hom_{TOP_M}((\mathcal{F}_{M'},\varphi_f),\underline{i}_*(L_{Z'},\phi_g))$$
de sorte que le foncteur $\underline{i}_*$ défini ci-dessus est bien l'image directe du morphisme $\underline{i}$. Le foncteur $\underline{i}_*$ est pleinement fidèle, car tous les $i'_*$ le sont. En d'autres termes, $\underline{i}:TOP_Z\rightarrow TOP_M$ est un plongement.

Montrons que l'image du plongement $\underline{i}$ est le sous-topos fermé complémentaire du sous-topos ouvert donné par $\underline{j}:TOP_U\rightarrow TOP_M $. Il faut montrer que l'image essentielle de $\underline{i}_*$ est précisément la sous-catégorie strictement pleine de $TOP_M$ formée des objets dont la restriction à $TOP_U$ est l'objet final. Soit $(L_{Z'},\phi_g)$ un objet de $TOP_Z$. Alors $\underline{j}^*\underline{i}_*(L_{Z'},\phi_g)=(F_{U'},\varphi_g)$, avec $F_{U'}=i'_*L_{\emptyset}$ où $i':\emptyset\rightarrow U'$. Donc $F_{U'}$ est l'objet final de $Sh(U')$ quel que soit l'espace $U'$ au-dessus de $U$, et $\underline{j}^*\underline{i}_*(L_{Z'},\phi_g)$ est l'objet final de $TOP_U$. Montrons maintenant que tout objet de $TOP_M$ dont la restriction à $TOP_U$ est l'objet final est dans l'image essentielle de $\underline{i}_*$. Soit $(F_{M'},\varphi_f)$ un objet de $TOP_M$ tel que $j^*(F_{M'},\varphi_f)$ est l'objet final de $TOP_U$. Pour tout $M'\rightarrow M$, on note $U'=M'\times_MU\rightarrow U$. La flèche $U'\rightarrow M'$ est une immersion ouverte, donc $F_{M'}\mid U'\simeq F_{U'}$ est l'objet final de $Sh(U')$. Il existe donc (essentiellement) un unique faisceau $L_{Z'}$ sur $Z'$ tel que $i'_*L_{Z'}=F_{M'}$, car l'image essentielle de $i'_*:Sh(Z')\rightarrow Sh(M')$ est la sous-catégorie strictement pleine de $Sh(M')$ formée des objets de $Sh(M')$ dont la restriction à $U'$ est l'objet final. Pour une fonction continue
$f:M''\rightarrow M'$ au-dessus de $M$, le morphisme $\varphi_f:f^*F_{M'}\rightarrow F_{M''}$ induit un morphisme
$f^*i'_*L_{Z'}\rightarrow i''_*L_{Z''}$.
\`A nouveau par changement de base, on obtient
$$i''_*g^*L_{Z'}\simeq f^*i'_*L_{Z'}\rightarrow i''_*L_{Z''}$$
car $i'$ est une immersion fermée. Comme le foncteur $i''_*$ est pleinement fidèle, ce morphisme est induit par un unique morphisme
$$\phi_g:g^*L_{Z'}\rightarrow L_{Z''}$$
qui est un isomorphisme lorsque $g$ est un étalement.
On obtient un objet $(L_{Z'},\phi_g)$ de $TOP_Z$ tel que $\underline{i}_*(L_{Z'},\phi_g)\simeq (F_{M'},\varphi_f)$. Plus précisément, on a $(F_{M'},\varphi_f)\simeq \underline{i}_*\underline{i}^*(F_{M'},\varphi_f)$. On a montré que \emph{$\underline{i}$ est un plongement fermé}.

L'équivalence $\mathcal{S}(M)\simeq TOP(M)$ est donnée par $$\mathcal{S}(M)=\mathcal{T}/yM
=\widetilde{(Top,\mathcal{J}_{ouv})}/yM\simeq \widetilde{(Top/M,\mathcal{J}_{ouv})}=TOP(M).$$
L'équivalence $TOP(M)\simeq TOP_M$ est définie comme suit. Si $\mathcal{F}$ est un faisceau sur $(Top/M,\mathcal{J}_{op})$ alors pour tout espace $M'$ au-dessus de $M$, $\mathcal{F}$ définit, par restriction au site $Ouv(M')$ des ouverts de $M'$, un petit faisceau $\mathcal{F}_{M'}$ sur $M'$. Pour une flèche $f:M''\rightarrow M'$, on a un morphisme $\mathcal{F}_{M'}\rightarrow f_*\mathcal{F}_{M''}$. Ce dernier est donné par $$\mathcal{F}(V')=\mathcal{F}_{M'}(V')\longrightarrow
f_*\mathcal{F}_{M''}(V')=\mathcal{F}_{M''}(M''\times_{M'}V')=\mathcal{F}(M''\times_{M'}V')$$
pour tout ouvert $V'$ de $M'$. On obtient alors $\varphi_f$ par adjonction.
Le résultat est maintenant donné par la commutativité du diagramme suivant :
\[ \xymatrix{
\mathcal{S}(U)\ar[d]^{\simeq}\ar[r] &\mathcal{S}(M)\ar[d]^{\simeq} &\mathcal{S}(Z)\ar[d]^{\simeq}\ar[l]\\
TOP(U)\ar[d]^{\simeq}\ar[r] &TOP(M)\ar[d]_{\alpha_M}^{\simeq} &TOP(Z)\ar[d]_{\alpha_Z}^{\simeq}\ar[l]_{TOP(i)}\\
TOP_U\ar[r] &TOP_M &TOP_Z\ar[l]_{\underline{i}}
} \]
On vérifie la commutativité du carré en bas à droite. Le reste sera laissé au lecteur. Soit $\mathcal{L}$ un objet de $TOP(Z)$. On pose
$$\alpha_{M*}\circ TOP(i)_*\mathcal{L}=(F_{M'},\varphi_f).$$
$F_{M'}$ est la restriction de $TOP(i)_*\mathcal{L}$ au site $Ouv(M')$. Pour un ouvert $V'$ de $M'$, on a
$$F_{M'}(V'):=TOP(i)_*\mathcal{L}(V')=\mathcal{L}(Z\times_MV').$$
Pour une flèche $f:M''\rightarrow M'$, le morphisme $\varphi_f$ est donné par la flèche
\begin{equation}\label{une map}
\mathcal{L}(V'\times_{M}Z)=TOP(i)_*\mathcal{F}(V')\rightarrow TOP(i)_*\mathcal{L}(M''\times_{M'}V')=\mathcal{L}(M''\times_{M'}V'\times_{M}Z)
\end{equation}
définie pour tout ouvert $V'$ de $M'$.

On considère maintenant $\underline{i}_*\circ \alpha_{Z*}\mathcal{L}$. On pose $$\underline{i}_*\circ \alpha_{Z*}\mathcal{L}=(F'_{M'},\varphi'_f).$$
On a $\alpha_{Z*}\mathcal{L}=(L_{Z'},\phi_g)$, où $L_{Z'}$ est la restriction de $\mathcal{L}$ au site $Ouv(Z')$. Alors $\underline{i}_*(L_{Z'},\phi_g)=(F'_{M'},\varphi'_f)$ avec $F'_{M'}=i'_*L_{Z'}$ pour $Z':=M'\times_MZ$. Pour un ouvert $V'$ de $M'$, on a $$F'_{M'}(V')=i'_*L_{Z'}(V')=L_{Z'}(Z'\times_{M'}V')=\mathcal{L}(Z'\times_{M'}V')=\mathcal{L}(Z\times_MV').$$
Pour une flèche $f:M''\rightarrow M'$, le morphisme $\varphi'_f:f^*F'_{M'}\rightarrow F'_{M''}$ correspond à $$F'_{M'}=i'_*L_{Z'}\longrightarrow f_*F'_{M''}=f_*i''_*L_{Z''}$$
qui est donné par la flèche
\begin{equation}\label{deux map}
i'_*L_{Z'}(V')=\mathcal{L}(V'\times_{M'}Z')\longrightarrow f_*i''_*L_{Z''}(V')=\mathcal{L}(V'\times_{M'}Z'')
\end{equation}
définie pour tout ouvert $V'$ de $M'$. Les flèches (\ref{une map}) et (\ref{deux map}) sont les mêmes, car elles sont données par les morphismes de restriction du faisceau $\mathcal{L}$. On en déduit un isomorphisme de foncteurs $\alpha_{M*}\circ TOP(i)_*\simeq \underline{i}_*\circ \alpha_{Z*}$. Donc le carré en bas à droite dans le diagramme précédent est bien commutatif. La commutativité de ce diagramme, le fait que $\underline{i}$ soit le plongement fermé complémentaire de $\underline{j}$ et le fait que les morphismes verticaux soient des équivalences montrent que l'on a une décomposition ouverte-fermée :
$$\mathcal{S}(U)\longrightarrow \mathcal{S}(M) \longleftarrow \mathcal{S}(Z)$$

\begin{lem}\label{plongementferme-equi} Soient $G$ un groupe topologique opérant continûment sur un espace $M$, $i:Z\hookrightarrow M$ un sous-espace fermé stable sous l'action de $G$ et $U$ l'ouvert complémentaire . Alors le morphisme de topos $\mathcal{S}(G,Z)\rightarrow\mathcal{S}(G,M)$ est le plongement fermé complémentaire du sous-topos ouvert $\mathcal{S}(G,U)\rightarrow\mathcal{S}(G,M)$.
\end{lem}
L'idée de la preuve est la suivante. On considère les pull-backs suivants, obtenus par localisation :
\[ \xymatrix{
\mathcal{S}(Z)\ar[d]\ar[r] &\mathcal{S}(M)\ar[r]\ar[d] &\mathcal{T}\ar[d]\\
\mathcal{S}(G,Z)\ar[r] &\mathcal{S}(G,M)\ar[r] & B_G
} \]
Le résultat est connu "en haut" par le lemme précédent. On va donc se ramener à cette situation non équivariante par descente. On a un morphisme du topos simplicial
\begin{equation}\label{top-simplicial1}
\mathcal{S}(G\times G\times Z)\rightrightarrows\rightarrow \mathcal{S}(G\times Z)\rightrightarrows\leftarrow\mathcal{S}(Z)
\end{equation}
dans le topos simplicial
\begin{equation}\label{top-simplicial2}
\mathcal{S}(G\times G\times M)\rightrightarrows\rightarrow \mathcal{S}(G\times M)\rightrightarrows\leftarrow\mathcal{S}(M)
\end{equation}
D'après le lemme précédent, ce morphisme de topos simpliciaux est composé de trois plongements fermés, qui sont les plongements fermés complémentaires des sous-topos ouverts donnés par le morphisme du topos simplicial
\begin{equation}\label{top-simplicial3}
\mathcal{S}(G\times G\times U)\rightrightarrows\rightarrow \mathcal{S}(G\times U)\rightrightarrows\leftarrow\mathcal{S}(U)
\end{equation}
dans le topos simplicial (\ref{top-simplicial2}). On note $Desc(G,Z)$, $Desc(G,M)$ et $Desc(G,U)$ les topos de descente des topos simpliciaux $(\ref{top-simplicial1})$, $(\ref{top-simplicial2})$ et $(\ref{top-simplicial3})$ respectivement. On note $d_{M,\bullet}$ les morphismes de structure du topos simplicial (\ref{top-simplicial2}). Rappelons que $Desc(G,M)$ est la catégorie des objets de $\mathcal{F}$ de $\mathcal{S}(M)$ munis d'une donnée de descente, i.e. d'un isomorphisme $d^*_{M,0}F\rightarrow d_{M,1}^*F$ satisfaisant les conditions habituelles. On note enfin $i_0:\mathcal{S}(Z)\rightarrow \mathcal{S}(M)$ et $i_1:\mathcal{S}(G\times Z)\rightarrow \mathcal{S}(G\times M)$. Les topos de descente $Desc(G,Z)$, $Desc(G,M)$ et $Desc(G,U)$ sont canoniquement équivalents à $\mathcal{S}(G,Z)$,  $\mathcal{S}(G,M)$ et $\mathcal{S}(G,U)$ respectivement.

L'image directe du morphisme $I:Desc(G,Z)\rightarrow Desc(G,M)$ peut se définir comme suit. Soit $(\mathcal{L},b:d_{Z,0}^*\mathcal{L}\rightarrow d_{Z,1}^*\mathcal{L})$ un objet de $Desc(G,Z)$. On considère $i_{0,*}\mathcal{L}$. On a un isomorphisme $i_{1,*}d_{Z,0}^*\mathcal{L}\rightarrow i_{1,*}d_{Z,1}^*\mathcal{L}$. Mais les morphismes de topos $d_{M,0}$ et $d_{M,1}$ sont des morphismes de localisation. On en déduit des isomorphismes $d_{M,0}^*i_{0,*}\simeq i_{1,*}d_{Z,0}^*$ et $d_{M,1}^*i_{0,*}\simeq i_{1,*}d_{Z,1}^*$. On a donc un isomorphisme
$$a:d_{M,0}^*i_{0,*}\mathcal{L} \simeq i_{1,*}d_{Z,0}^*\mathcal{L}\rightarrow i_{1,*}d_{Z,1}^*\mathcal{L}\simeq d_{M,1}^*i_{0,*}\mathcal{L}.$$
On obtient un objet $(i_{0,*}\mathcal{L},a)$ de $Desc(G,M)$ qui est l'image directe de $(\mathcal{L},b)$. Il est clair que $I$ est un plongement (i.e. $I_*$ est pleinement fidèle), et que l'image du morphisme $I$ est contenu dans le fermé complémentaire du plongement ouvert $Desc(G,U)\rightarrow Desc(G,M)$. Il reste à montrer que l'image de $I$ est précisément ce sous-topos fermé.

Soit $(\mathcal{F},a:d^*_{M,0}\mathcal{F}\rightarrow d_{M,1}^*\mathcal{F})$ un objet de $Desc(G,M)$ dont l'image inverse dans $Desc(G,U)$ est l'objet final. En particulier $\mathcal{F}\mid U$ est l'objet final de $\mathcal{S}(U)$. D'après le lemme précédent, $\mathcal{F}=i_{0,*}\mathcal{L}$ pour un objet $\mathcal{L}$ de $\mathcal{S}(Z)$. On a donc une flèche $d_{M,0}^*i_{0,*}\mathcal{L}\rightarrow d_{M,1}^*i_{0,*}\mathcal{L}$. A nouveau, les morphismes $d_{M,0}$ et $d_{M,1}$ sont des morphismes de localisation. On en déduit $d_{M,0}^*i_{0,*}\simeq i_{1,*}d_{Z,0}^*$, $d_{M,1}^*i_{0,*}\simeq i_{1,*}d_{Z,1}^*$ et un isomorphisme
$$i_{1,*}d_{Z,0}^*\mathcal{L}\simeq d_{M,0}^*i_{0,*}\mathcal{L}\rightarrow d_{M,1}^*i_{0,*}\mathcal{L}\simeq i_{1,*}d_{Z,1}^*\mathcal{L}$$
Comme $i_{1,*}$ est pleinement fidèle, on obtient un objet
$(\mathcal{L},b:d_{Z,0}^*\mathcal{L}\rightarrow d_{Z,1}^*\mathcal{L})$  de $Desc(G,Z)$ dont l'image directe dans $Desc(G,M)$ est $(\mathcal{F},a:d^*_{M,0}\mathcal{F}\rightarrow d_{M,1}^*\mathcal{F})$.

On considère maintenant le diagrame commutatif suivant.
\[ \xymatrix{
Desc(G,U)\ar[d]^{\simeq}\ar[r] &Desc(G,M)\ar[d]^{\simeq} &Desc(G,Z)\ar[d]^{\simeq}\ar[l]\\
\mathcal{S}(G,U)\ar[r] &\mathcal{S}(G,M) &\mathcal{S}(G,U)\ar[l]
} \]
On vient de montrer que la ligne du haut est une décomposition ouverte-fermée. Le lemme \ref{plongementferme-equi} suit car les flèches verticales sont des équivalences. Ceci achève la preuve du théorème \ref{orbiteferme+ptsfixes}.

\end{f-proof}

\subsection{Le morphisme $\mathcal{S}(\mathbb{R},\mathrm{M}_{\bar{X}})\rightarrow\bar{X}_{et}$}
Le site étale d'Artin-Verdier $Et_{\bar{X}}$ est défini dans \cite{Artin-Verdier}. Le topos étale d'Artin-Verdier $\bar{X}_{et}$
est la catégorie des faisceaux d'ensembles sur le site $Et_{\bar{X}}$. Supposons que le foncteur conjecturé par C. Deninger existe. Par
restriction, on obtient un foncteur
$$
\fonc{\gamma^*}{Et_{\bar{X}}}{{Top}^{\mathbb{R}}/_{\mathrm{M}_{\bar{X}}}}
{\bar{U}\rightarrow\bar{X}}{\mathrm{M}_{\bar{U}}\rightarrow\mathrm{M}_{\bar{X}}}
$$
\begin{hypothese}\label{hyp2} Nous supposons que les conditions suivantes sont satisfaites.
\begin{itemize}
\item Le foncteur $\gamma^*$
existe et commute aux produits fibrés.
\item Un morphisme étale $\bar{V}\rightarrow\bar{U}$ est envoyé sur un homémorphisme local.
\item Une famille surjective de morphismes étales de schémas est envoyée sur une famille surjective d'applications continues.
\end{itemize}
\end{hypothese}
Une famille surjective de morphismes étales de schémas
$\{\bar{U_i}\rightarrow\bar{U};\,\,i\in I\}$ induit une
famille surjective d'étalements $\mathbb{R}$-équivariants
$\{\mathrm{M}_{\bar{U_i}}\rightarrow\mathrm{M}_{\bar{U}};\,\,i\in
I\}$. Puisqu'un étalement admet des sections locales au-dessus de
son image, le foncteur $\gamma^*$ envoie un recouvrement pour la
prétopologie étale sur un recouvrement pour la prétopologie des
sections locales. Ainsi, on a un morphisme de sites exacts à gauche
$$\gamma^*:(Et_{\bar{X}};\mathcal{J}_{et})\longrightarrow
({Top}^{\mathbb{R}}/_{(\mathbb{R},\mathrm{M}_{\bar{X}})};\mathcal{J}_{ls}).$$
La proposition suivante en découle immédiatement.
\begin{prop} \label{mapgammapourM}
Si l'hypothèse \ref{hyp2} est satisfaite, alors on a un morphisme de topos
$$\gamma:\mathcal{S}(\mathbb{R};\mathrm{M}_{\bar{X}})\longrightarrow\bar{X}_{et}.$$
\end{prop}

Soit $v$ un point fermé de $\bar{X}$ ($v$ est une place ultramétrique ou archimédienne du corps de nombres $K$). Lorsque $v$ est ultramétrique, on note $W_{k(v)}$ le sous-groupe de $\mathbb{R}$ engendré par $log(N(v))$. Lorsque $v$ est archimédienne, on note $W_{k(v)}=\mathbb{R}$ et $G_{k(v)}=1$. On a un morphisme canonique $W_{k(v)}\rightarrow G_{k(v)}$ induisant un morphisme $\alpha_v:B_{W_{k(v)}}\rightarrow B^{sm}_{G_{k(v)}}$. On note aussi
$$u_v:B^{sm}_{k(v)}\simeq Spec(k(v))_{et}\longrightarrow\bar{X}_{et}$$
le plongement fermé induit par l'inclusion fermée de schémas
$Spec(k(v))\rightarrow\bar{X}$, pour tout point fermé $v$ de $\bar{X}$.

\begin{prop}\label{commdiagpourM}
Soit $v$ une place ultramétrique ou archimédienne du corps de nombres $K$. Si l'hypothèse \ref{hyp2}
est satisfaite, alors on a un diagramme commutatif de topos
\[ \xymatrix{
 B_{W_{k(v)}} \ar[d]_{i_v} \ar[r]^{\alpha_v} &B^{sm}_{G_{k(v)}} \ar[d]_{u_v}   \\
 \mathcal{S}(\mathbb{R};\mathrm{M}_{\bar{X}})\ar[r]^{\gamma}&\bar{X}_{et}
} \]
De plus, le morphisme $i_v$ est un plongement fermé.
\end{prop}
\begin{f-proof}
On traite d'abord le cas d'une place ultramétrique $v$.
Le morphisme de topos $\alpha_v:B_{W_{k(v)}}\rightarrow B^{sm}_{G_{k(v)}}$,
défini dans la proposition \ref{morphisme alpha}, est induit par le morphisme de sites exacts à gauche
$$
\fonc{\alpha^*_v}{(Et_{Spec(k(v))};\mathcal{J}_{et})}{(Top^{\mathbb{R}}/_{\mathbb{M}_{N(v)}};\mathcal{J}_{ls})}
{\mathbb{F}_q/k(v)}{\mathbb{M}_q/\mathbb{M}_{N(v)}}
$$
Le plongement fermé $u_v:B^{sm}_{k(v)}\rightarrow\bar{X}_{et}$,
induit par l'inclusion fermée de schémas
$Spec(k(v))\rightarrow\bar{X}$, est défini par le morphisme de sites exacts à gauche
$$
\fonc{u_v^*}{(Et_{\bar{X}};\mathcal{J}_{et})}{(Et_{Spec(k(v))};\mathcal{J}_{et})}{\bar{U}}{\bar{U}\times_{\bar{X}}Spec(k(v))}
$$
D'autre part, le morphisme de topos
$i_v:B_{W_{k(v)}}\rightarrow\mathcal{S}(\mathbb{R};\mathrm{M}_{\bar{X}})$
est induit par le morphisme de sites exacts à gauche
$$
\fonc{i_v^*}{(Top^{\mathbb{R}}/_{\mathrm{M}_{\bar{X}}};\mathcal{J}_{ls})}{(Top^{\mathbb{R}}/_{\mathbb{M}_{N(v)}};\mathcal{J}_{ls})}
{Z}{Z\times_{\mathrm{M}_{\bar{X}}}\mathbb{M}_{N(v)}}
$$
où $\mathbb{R}$ opère diagonalement sur
$Z\times_{\mathrm{M}_{\bar{X}}}\mathbb{M}_{N(v)}$. Alors, le
diagramme suivant de sites est commutatif (car le foncteur $\mathfrak{D}$ est supposé commuter aux produits fibrés).
\[ \xymatrix{
{Top}^{\mathbb{R}}/_{\mathbb{M}_{N(v)}}   &\ar[l]_{\alpha_v^*} Et_{Spec(k(v))}   \\
 {Top}^{\mathbb{R}}/_{(\mathbb{R},\mathrm{M}_{\bar{X}})}\ar[u]_{i_v^*}
 & \ar[l]_{\gamma^*}\ar[u]_{u_v^*}Et_{\bar{X}}
} \]
Le diagramme de topos correspondant est donc lui aussi commutatif.

Soit maintenant $v=\mathfrak{p}\in X_{\infty}\subset\bar{X}$ une
place archimédienne de $K$, et $m_{\mathfrak{p}}\in
\mathrm{M}_{\bar{X}}$ le point fixe correspondant. Le morphisme de topos
$$\alpha_{\mathfrak{p}}:\mathcal{S}(\mathbb{R};\{m_{\mathfrak{p}}\})=B_{\mathbb{R}}\longrightarrow B^{sm}_{G_{k(\mathfrak{p})}}=\underline{Set}$$
est induit par le morphisme de sites
exacts à gauche
$\alpha^*_{\mathfrak{p}}:(\underline{Set};\mathcal{J}_{can})\rightarrow(Top^{\mathbb{R}};\mathcal{J}_{ls})$
qui envoie un ensemble $E$ sur l'espace topologique discret $E$ sur
lequel $\mathbb{R}$ opère trivialement. Le morphisme
$u_{\mathfrak{p}}:\underline{Set}\rightarrow\bar{X}_{et}$,
induit par l'inclusion fermée
$\mathfrak{p}\rightarrow\bar{X}$, est donné par le morphisme de sites exacts à gauche
$$
\fonc{u_{\mathfrak{p}}^*}{(Et_{\bar{X}};\mathcal{J}_{et})}{(\underline{Set};\mathcal{J}_{can})}
{\bar{U}}{\bar{U}\times_{\bar{X}}\mathfrak{p}}
$$
Le morphisme de topos $i_{\mathfrak{p}}:B_{\mathbb{R}}\rightarrow\mathcal{S}(\mathbb{R};\mathrm{M}_{\bar{X}})$
est induit par le morphisme de sites exacts à gauche
$$
\fonc{i_{\mathfrak{p}}^*}{(Top^{\mathbb{R}}/_{\mathrm{M}_{\bar{X}}};\mathcal{J}_{ls})}
{(Top^{\mathbb{R}};\mathcal{J}_{ls})}
{Z}{Z\times_{\mathrm{M}_{\bar{X}}}m_\mathfrak{p}}
$$
On obtient à nouveau un diagramme commutatif de sites et le résultat suit. Le fait que $i_v$ est un plongement fermé est donné par le théorème \ref{orbiteferme+ptsfixes}
\end{f-proof}

\section{Propriétés analogues du topos
Weil-étale}

Soit $K$ un corps de nombres, $X=Spec(\mathcal{O}_K)$, $X_{\infty}$ l'ensemble des places archimédiennes de $K$, et $\bar{X}=(X,X_{\infty})$. Le topos Weil-étale $\bar{X}_W$ est défini dans \cite{Flach-moi} comme une modification de la définition donnée par Lichtenbaum dans \cite{Lichtenbaum}. Les résultats suivants sont démontrés dans \cite{Flach-moi}.

\begin{thm}\label{thm-weiletaletopos}
On a les résultats suivants.
\begin{enumerate}
\item On a un morphisme canonique $\mathfrak{f}:\bar{X}_W\rightarrow B_{\mathbb{R}}$.
\item Si $v$ est un point fermé de $\bar{X}$ (i.e. une place ultramétrique ou archimédienne de $K$), alors on a un plongement fermé $i_v:B_{W_{k(v)}}\rightarrow\bar{X}_W$.
De plus, la composition $\mathfrak{f}\circ
i_v:B_{W_{k(v)}}\rightarrow B_{\mathbb{R}}$ est induite par le morphisme canonique $l_v:W_{k(v)}\rightarrow\mathbb{R}$.
\item On a un morphisme canonique
$\gamma:\bar{X}_W\rightarrow \bar{X}_{et}$.
\item Si $v$ est un point fermé de $\bar{X}$, on a un diagramme commutatif
\[ \xymatrix{
 B_{W_{k(v)}} \ar[d]_{i_v} \ar[r]^{\alpha_v} &B^{sm}_{G_{k(v)}} \ar[d]_{u_v}   \\
 \bar{X}_W\ar[r]^{\gamma}&\bar{X}_{et}
} \]
De plus, ce diagramme est un pull-back.
\end{enumerate}
\end{thm}
\begin{rem}En supposant que le système dynamique $\mathrm{M}_{\bar{X}}$ existe,
les propriétés précédentes satisfaites par le topos Weil-étale sont aussi satisfaites par le topos $\mathcal{S}(\mathbb{R},\mathrm{M}_{\bar{X}})$. En effet, $(1)$ est donnée par la définition \ref{def-flot}, $(2)$ est donnée par le théorème \ref{orbiteferme+ptsfixes}, $(3)$ est donnée par la propriété \ref{mapgammapourM} et la commutativité du diagramme $(4)$ est donnée par la propriété \ref{commdiagpourM}. On interprète donc le morphisme $\mathfrak{f}:\bar{X}_W\rightarrow B_{\mathbb{R}}$ comme un flot. Lorsque $v$ est une place archimédienne de $K$, la flèche $i_v:B_{W_{k(v)}}\rightarrow\bar{X}_W$ est vue comme une orbite fermée de longueur $log(N(v))$. Lorsque $v$ est une place archimédienne de $K$, la flèche $i_v:B_{W_{k(v)}}\rightarrow\bar{X}_W$ est vue comme un point fixe du flot.
\end{rem}

\section{Schémas arithmétiques de dimension supérieure}

Soit $\mathcal{X}\rightarrow Spec(\mathbb{Z})$ un schéma arithmétique régulier, connexe, plat et propre sur $Spec(\mathbb{Z})$. On définit $\mathcal{X}_{\infty}:=\mathcal{X}(\mathbb{C})/G_{\mathbb{R}}$, où $\mathcal{X}(\mathbb{C})$ est muni de la topologie complexe et $\mathcal{X}_{\infty}$ de la topologie quotient. On considère $\bar{\mathcal{X}}:=(\mathcal{X},\mathcal{X}_{\infty})$. Dans cette section, on met en évidence certaines propriétés partagées par le topos Weil-étale $\bar{\mathcal{X}}_W$ défini dans \cite{Flach-moi} et le système dynamique conjecturalement associé à $\bar{\mathcal{X}}$.

\subsection{Le topos Weil-étale $\bar{\mathcal{X}}_W$}
Le site étale d'Artin-Verdier $Et_{\bar{\mathcal{X}}}$ est défini dans \cite{Artin-Verdier}. Le topos étale d'Artin-Verdier $\bar{\mathcal{X}}_{et}$
est la catégorie des faisceaux d'ensembles sur le site $Et_{\bar{\mathcal{X}}}$. Ce topos est étudié dans \cite{Flach-moi}. On a un morphisme $\bar{\mathcal{X}}_{et}\rightarrow\overline{Spec(\mathbb{Z})}_{et}$. On définit alors le topos Weil-étale de $\bar{\mathcal{X}}$ comme le produit fibré
$$\bar{\mathcal{X}}_W:=\bar{\mathcal{X}}_{et}\times_{\overline{Spec(\mathbb{Z})}_{et}}\overline{Spec(\mathbb{Z})}_{W}$$
dans la 2-catégorie des topos, où $\overline{Spec(\mathbb{Z})}_{W}$ est le topos Weil-étale de $\overline{Spec(\mathbb{Z})}$ utilisé dans la section précédente. La deuxième projection
$\bar{\mathcal{X}}_W\rightarrow \overline{Spec(\mathbb{Z})}_{W}$
induit un morphisme
$$\mathfrak{f}_{\bar{\mathcal{X}}}:\bar{\mathcal{X}}_W\rightarrow \overline{Spec(\mathbb{Z})}_{W}\rightarrow B_{\mathbb{R}}.$$
Un point fermé $x\in\mathcal{X}$ fournit à nouveau un plongement fermé
$$i_x:B_{W_{k(x)}}=\mathcal{S}(\mathbb{R},\mathbb{R}/log(N(x))\mathbb{Z})\longrightarrow\bar{\mathcal{X}}_W$$
de sorte que la composition $\mathfrak{f}_{\bar{\mathcal{X}}}\circ i_x:B_{W_{k(x)}}\rightarrow  B_{\mathbb{R}}$ soit induite par le morphisme canonique $W_{k(x)}\rightarrow\mathbb{R}$.

\subsection{Le topos  $\mathcal{S}(\mathbb{R},\mathrm{M}_{\bar{\mathcal{X}}})$}
Supposons que l'on puisse associer au schéma arithmétique $\bar{\mathcal{X}}$ un système dynamique $\mathrm{M}_{\bar{\mathcal{X}}}$, de la forme décrite dans (\cite{Deninger3} Dictionary 4 , part 2). Le morphisme $\bar{\mathcal{X}}\rightarrow\overline{Spec(\mathbb{Z})}$ induit en particulier une application continue $\mathbb{R}$-équivariante $\mathrm{M}_{\bar{\mathcal{X}}}\rightarrow\mathrm{M}_{\overline{Spec(\mathbb{Z})}}$. On en déduit un morphisme
$$\mathfrak{f}_{\mathrm{M}_{\bar{\mathcal{X}}}}:\mathcal{S}(\mathbb{R},\mathrm{M}_{\bar{\mathcal{X}}})\rightarrow \mathcal{S}(\mathbb{R},\mathrm{M}_{\overline{Spec(\mathbb{Z})}})\rightarrow B_{\mathbb{R}}.$$
Une orbite fermée $\gamma_x$ correspondant à un point fermé $x\in\mathcal{X}$ fournit un plongement fermé
$$i_x:B_{W_{k(x)}}=\mathcal{S}(\mathbb{R},\mathbb{R}/log(N(x))\mathbb{Z})\longrightarrow \mathcal{S}(\mathbb{R},\mathrm{M}_{\bar{\mathcal{X}}})$$
de sorte que la composition $\mathfrak{f}_{\mathrm{M}_{\bar{\mathcal{X}}}}\circ i_x:B_{W_{k(x)}}\rightarrow  B_{\mathbb{R}}$ soit induite par le morphisme canonique $W_{k(x)}\rightarrow\mathbb{R}$.

\subsection{Structure du topos $\bar{\mathcal{X}}_{W}$ au-dessus des points fermés de $Spec(\mathbb{Z})$}

Soit $p$ un point fermé de $Spec(\mathbb{Z})$. Il lui correspond un plongement fermé
$$i_v:B_{W_{\mathbb{F}_p}}=\mathcal{S}(\mathbb{R},\mathbb{R}/log(p)\mathbb{Z})\longrightarrow\overline{Spec(\mathbb{Z})}_W.$$
La fibre de $\bar{\mathcal{X}}_{W}$ au-dessus du point fermé $p\in Spec(\mathbb{Z})$ est définie comme l'image inverse du sous-topos fermé $Im(i_v)$. Elle s'identifie au produit fibré
$$\bar{\mathcal{X}}_{W}\times_{\overline{Spec(\mathbb{Z})}_W}B_{W_{\mathbb{F}_p}}.$$
Alors on a des équivalences canoniques (cf. \cite{Flach-moi})
\begin{equation}\label{fibredeWETenp}
\bar{\mathcal{X}}_{W}\times_{\overline{Spec(\mathbb{Z})}_W}B_{W_{\mathbb{F}_p}}\simeq (\mathcal{X}\otimes\mathbb{F}_p)_{et}\times_{B^{sm}_{G_{\mathbb{F}_p}}}B_{W_{\mathbb{F}_p}}
\simeq(\mathcal{X}\otimes\mathbb{F}_p)_W.
\end{equation}
Ainsi la fibre du topos $\bar{\mathcal{X}}_{W}$ au-dessus du point fermé $p$ est équivalente au topos Weil-étale du schéma $\mathcal{X}\otimes\mathbb{F}_p$.

\subsection{Structure du topos $\mathcal{S}(\mathbb{R},\bar{\mathcal{X}})$ au-dessus  des points fermés de $Spec(\mathbb{Z})$}

Soit $p$ un point fermé de $Spec(\mathbb{Z})$. Il lui correspond une orbite fermée
$\gamma_p:\mathbb{R}/log(p)\mathbb{Z} \rightarrow\mathrm{M}_{\overline{Spec(\mathbb{Z})}}$
qui induit un plongement fermé de topos
$$i_v:B_{W_{\mathbb{F}_p}}=\mathcal{S}(\mathbb{R},\mathbb{R}/log(p)\mathbb{Z})\longrightarrow
\mathcal{S}(\mathbb{R},\mathrm{M}_{\overline{Spec(\mathbb{Z})}}).$$
La fibre de l'espace $\mathrm{M}_{\bar{\mathcal{X}}}$ au-dessus du point fermé $p\in Spec(\mathbb{Z})$ (i.e. de l'orbite fermée $\gamma_p$) est naturellement définie comme le produit fibré (d'espaces topologiques)
$$
\mathrm{M}_{\bar{\mathcal{X}}}\times_{\mathrm{M}_{\overline{Spec(\mathbb{Z})}}}(\mathbb{R}/log(p)\mathbb{Z})
=\mathrm{M}_{\bar{\mathcal{X}}}\times_{\mathrm{M}_{\overline{Spec(\mathbb{Z})}}}\mathrm{M}_{\mathbb{F}_p}$$
D'après (\cite{SGA4} IV Proposition 5.11) on a une équivalence de topos
$$\mathcal{S}(\mathbb{R},\mathrm{M}_{\bar{\mathcal{X}}}\times_{\mathrm{M}_{\overline{Spec(\mathbb{Z})}}}
(\mathbb{R}/log(p)\mathbb{Z}))\simeq
\mathcal{S}(\mathbb{R},\mathrm{M}_{\bar{\mathcal{X}}})\times_{\mathcal{S}(\mathbb{R},\mathrm{M}_{\overline{Spec(\mathbb{Z})}})}
B_{W_{\mathbb{F}_p}}$$
On suppose à nouveau que le foncteur envisagé par Deninger existe, disons sur la catégorie des schémas de type fini séparés sur $Spec(\mathbb{Z})$ (cf. Section \ref{sectionmorphisme}),  et qu'il commute aux produits fibrés. Alors on a
$$\mathrm{M}_{\mathcal{X}\otimes\mathbb{F}_p}=\mathrm{M}_{\bar{\mathcal{X}}}\times_{\mathrm{M}_{\overline{Spec(\mathbb{Z})}}}\mathrm{M}_{\mathbb{F}_p}.$$
On obtient donc une équivalence de topos
\begin{equation}\label{fibredeMenp}
\mathcal{S}(\mathbb{R},\mathrm{M}_{\mathcal{X}\otimes\mathbb{F}_p})\simeq
\mathcal{S}(\mathbb{R},\mathrm{M}_{\bar{\mathcal{X}}})\times_{\mathcal{S}(\mathbb{R},\mathrm{M}_{\overline{Spec(\mathbb{Z})}})}
B_{W_{\mathbb{F}_p}}.
\end{equation}

\begin{rem}
D'après le théorème \ref{propmorphisme}, on a de plus un morphisme du topos (\ref{fibredeMenp}), qui est la fibre de $\mathcal{S}(\mathbb{R},\mathrm{M}_{\bar{\mathcal{X}}})$ au-dessus de $p\in Spec(\mathbb{Z})$, dans le topos (\ref{fibredeWETenp}), qui est la fibre de $\bar{\mathcal{X}}_W$ au-dessus de $p\in Spec(\mathbb{Z})$. Ce morphisme
$$d:\mathcal{S}(\mathbb{R};\mathrm{M}_{\mathcal{X}\otimes\mathbb{F}_p})\longrightarrow (\mathcal{X}\otimes\mathbb{F}_p)_W$$
est compatible au flot, aux orbites fermées et au feuilletage (voir la section \ref{sectionmorphisme}).
\end{rem}

\subsection{Structure du topos $\bar{\mathcal{X}}_{W}$ au-dessus de la place archimédienne}
Soit $\infty\in\overline{Spec(\mathbb{Z})}$ la place archimédienne de $\mathbb{Q}$. Il lui correspond un plongement fermé
$$i_\infty:B_{W_{k(\infty)}}=B_{\mathbb{R}}\longrightarrow\overline{Spec(\mathbb{Z})}_W.$$
La fibre de $\bar{\mathcal{X}}_{W}$ au-dessus de $\infty$ est définie comme l'image inverse du sous-topos fermé $Im(i_\infty)$. Elle s'identifie au produit fibré
$$\bar{\mathcal{X}}_{W}\times_{\overline{Spec(\mathbb{Z})}_W}B_{\mathbb{R}}.$$
Alors on a une équivalence
\begin{equation}\label{fibredeWET-en-infty}
{\mathcal{X}}_{\infty,W}:=\bar{\mathcal{X}}_{W}\times_{\overline{Spec(\mathbb{Z})}_W}B_{\mathbb{R}}\simeq Sh(\mathcal{X}_{\infty})\times B_{\mathbb{R}}.
\end{equation}
où $Sh(\mathcal{X}_{\infty})$ est la catégorie des (petits) faisceaux (i.e. des espaces étalés) sur l'espace topologique $\mathcal{X}_{\infty}$.
Ainsi la fibre du topos $\bar{\mathcal{X}}_{W}$ au-dessus de  la place archimédienne $\infty$ est équivalente au topos Weil-étale
$Sh(\mathcal{X}_{\infty})\times B_{\mathbb{R}}$.

\subsection{Structure du topos $\mathcal{S}(\mathbb{R},\bar{\mathcal{X}})$ au-dessus de la place archimédienne}
Supposons que l'on puisse associer un système dynamique $\mathrm{M}_{\bar{\mathcal{X}}}$ au schéma arithmétique $\bar{\mathcal{X}}$ de la forme décrite dans \cite{Deninger3}. Le morphisme $\bar{\mathcal{X}}\rightarrow\overline{Spec(\mathbb{Z})}$ induit une application continue $\mathbb{R}$-équivariante $\mathrm{M}_{\bar{\mathcal{X}}}\rightarrow\mathrm{M}_{\overline{Spec(\mathbb{Z})}}$.
On considère le point fixe $m_{\infty}\in \mathrm{M}_{\overline{Spec(\mathbb{Z})}}$ correspondant à la place archimédienne $\infty$. Il induit un plongement fermé de topos
$$B_{\mathbb{R}}\longrightarrow
\mathcal{S}(\mathbb{R},\mathrm{M}_{\overline{Spec(\mathbb{Z})}}).$$
La fibre de l'espace $\mathrm{M}_{\bar{\mathcal{X}}}$ au-dessus de $\infty\in \overline{Spec(\mathbb{Z})}$ (i.e. du point fixe $m_{\infty}$) est le produit fibré
$\mathrm{M}_{\bar{\mathcal{X}}}\times_{\mathrm{M}_{\overline{Spec(\mathbb{Z})}}}m_{\infty}$. On a une équivalence de topos
$$\mathcal{S}(\mathbb{R},\mathrm{M}_{\bar{\mathcal{X}}}\times_{\mathrm{M}_{\overline{Spec(\mathbb{Z})}}}
m_{\infty})\simeq
\mathcal{S}(\mathbb{R},\mathrm{M}_{\bar{\mathcal{X}}})\times_{\mathcal{S}(\mathbb{R},\mathrm{M}_{\overline{Spec(\mathbb{Z})}})}
B_{\mathbb{R}}$$
D'après (\cite{Deninger3} Dictionary 4 , part 2), on indentifie les espaces
$\mathrm{M}_{\bar{\mathcal{X}}}\times_{\mathrm{M}_{\overline{Spec(\mathbb{Z})}}}
m_{\infty}=\mathcal{X}_{\infty}$
avec l'action triviale de $\mathbb{R}$. Alors on a  $$\mathcal{S}(\mathbb{R},\mathrm{M}_{\bar{\mathcal{X}}}\times_{\mathrm{M}_{\overline{Spec(\mathbb{Z})}}}
m_{\infty})=\mathcal{S}(\mathbb{R},\mathcal{X}_{\infty}):=B_{\mathbb{R}}/_{y(\mathbb{R},\mathcal{X}_{\infty})}.$$
On obtient une équivalence
$$\mathcal{S}(\mathbb{R},\mathrm{M}_{\bar{\mathcal{X}}})\times_{\mathcal{S}(\mathbb{R},\mathrm{M}_{\overline{Spec(\mathbb{Z})}})}
B_{\mathbb{R}}\simeq \mathcal{S}(\mathbb{R},\mathcal{X}_{\infty})$$
On a de plus un morphisme $$\mathcal{S}(\mathbb{R},\mathcal{X}_{\infty})\rightarrow \mathcal{S}(\mathcal{X}_{\infty})\simeq TOP(\mathcal{X}_{\infty})\rightarrow Sh(\mathcal{X}_{\infty})$$
et un morphisme de localisation
$$\mathcal{S}(\mathbb{R},\mathcal{X}_{\infty}):=B_{\mathbb{R}}/_{y(\mathbb{R},\mathcal{X}_{\infty})}\longrightarrow B_{\mathbb{R}}.$$
Par définition du produit, on obtient un morphisme
$$\mathcal{S}(\mathbb{R},\mathcal{X}_{\infty})\longrightarrow Sh(\mathcal{X}_{\infty})\times B_{\mathbb{R}}.$$
Ces deux topos sont naturellement associés à l'action triviale de $\mathbb{R}$ sur l'espace $\mathcal{X}_{\infty}$. Ils ont la même cohomologie à valeurs dans un faisceau constant discret, mais ne sont pas équivalents. Une meilleure définition du topos associé à $\mathrm{M}_{\bar{\mathcal{X}}}$ (on pourrait par exemple considérer la catégorie des faisceaux sur $\mathrm{M}_{\bar{\mathcal{X}}}$ à valeurs dans $\mathcal{T}$) devrait faire du morphisme (\ref{mapMtoWETatinfty}) ci-dessous une équivalence.
\begin{rem}
On a un morphisme canonique de topos
\begin{equation}\label{mapMtoWETatinfty}
\mathcal{S}(\mathbb{R},\mathcal{X}_{\infty})\longrightarrow Sh(\mathcal{X}_{\infty})\times B_{\mathbb{R}}\simeq{\mathcal{X}}_{\infty,W}
\end{equation}
de la fibre de $\mathcal{S}(\mathbb{R},\mathrm{M}_{\bar{\mathcal{X}}})$ au-dessus de $\infty$ dans la fibre de $\bar{\mathcal{X}}_W$ au-dessus de $\infty$.
\end{rem}

\section{Appendice}\label{appendix}
Nous montrons dans cette section comment le morphisme $\mathfrak{f}:\bar{X}_W\rightarrow B_{\mathbb{R}}$, donné par le théorème \ref{thm-weiletaletopos}, permet de définir une action de $\mathbb{R}$ sur le topos $\bar{X}_W\times_{B_{\mathbb{R}}}\mathcal{T}$.

\subsection{}La donnée d'un espace topologique sobre $Z$ équivaut à celle du topos $Sh(Z)$ (cf. \cite{SGA4} IV. 4.2). Nous donnons une version équivariante de ce résultat. Considérons un espace topologique sobre $Z$ sur
lequel un groupe discret $G$ opère par automorphismes à gauche. On note
$$r:G\times Z\rightarrow Z$$
cette action, $\check{r}:Z\times G\rightarrow Z$ l'action à droite déduite de $r$.
Soient $Sh(G;Z)$ et $B^{sm}_G:=G-\underline{Set}$ les topos des
petits $G$-faisceaux sur $Z$ et sur le point $\{*\}$ respectivement. Le topos $Sh(G;Z)$ est équivalent à la catégorie des espaces étalés $\tilde{Z}$ sur $Z$ munis d'une action de $G$ tel que la projection $\tilde{Z}\rightarrow Z$ est $G$-équivariante. Une application continue $G$-équivariante $u:Z\rightarrow Z'$ induit un morphisme de topos $$Sh(u):Sh(G;Z)\rightarrow Sh(G;Z').$$ Par exemple, l'unique application continue $Z\rightarrow\{*\}$ induit un morphisme de topos
$$f:Sh(G;Z)\longrightarrow B^{sm}_G.$$
\begin{thm}
Soit $Z$ un espace sobre sur lequel un groupe discret $G$ opère. Le topos $Sh(G;Z)$ muni du morphisme $f:Sh(G;Z)\rightarrow B^{sm}_G$ détermine l'action de $G$ sur $Z$ à homéomorphisme $G$-équivariant près. Plus précisément, soient $Z$ et $Z'$ deux espaces sobres munis d'une action de $G$ et soit $\alpha:Sh(G;Z)\rightarrow Sh(G;Z')$ un morphisme de topos au-dessus $B^{sm}_G$. Alors il existe une unique application continue
$G$-équivariante $a:Z\rightarrow Z'$ induisant $\alpha$. De plus, $\alpha$ est une équivalence si et seulement si $a$ est un homéomorphisme.
\end{thm}
\begin{f-proof}
Si $E$ est un $G$-ensemble, alors $f^*E$ est la projection $E\times Z\rightarrow Z$, où $G$ opère diagonalement sur $E\times Z$.
Soit $EG$ l'objet de $B^{sm}_G$ donné par $G$ sur lequel $G$ opère par multiplication à gauche. Alors on a une équivalence canonique
$$Sh(G;Z)/_{f^*EG}\simeq Sh(Z)$$
Ce topos détermine l'espace $Z$  à homéomorphisme près. On peut aussi retrouver l'action de $G$ sur $Z$. Le groupe $G$ opère (par multiplication à droite) sur $EG$ dans $B^{sm}_G$. Comme $f^*$ commute aux limites projectives finies, le groupe $G$ opère à droite sur $f^*EG$ dans $Sh(G;Z)/_{f^*EG}\simeq Sh(Z)$. En particulier, un élément $g\in G$ définit un isomorphisme $f^*EG\rightarrow f^*EG$ qui, par transitivité des topos induits (cf. \cite{SGA4} IV.5.5), induit à son tour une équivalence
$$\tilde{g}: Sh(G;Z)/_{f^*EG}\longrightarrow Sh(G;Z)/_{f^*EG}.$$
On obtient une action à droite de $G$ sur le topos $Sh(Z)\simeq Sh(G;Z)/_{f^*EG}$, donnée par le morphisme
\begin{equation}\label{actionGsurSh}
\rho:Sh(Z)\times G:=\coprod_{G}Sh(Z)\longrightarrow Sh(Z)
\end{equation}
Cette action à droite de $G$ sur le topos $Sh(Z)$ est l'action induite par $r$ : on a
$$\rho=Sh(\check{r}).$$
En effet, quel que soit $g\in G$, on vérifie facilement que le diagramme de topos
\[ \xymatrix{
Sh(G;Z)/_{f^*EG}\ar[d]^{\widetilde{g}} \ar[r]^{\,\,\,\,\,\,\,\,\,\simeq} &Sh(Z) \ar[d]^{Sh(g^{-1})}   \\
Sh(G;Z)/_{f^*EG} \ar[r]^{\,\,\,\,\,\,\,\,\,\simeq}&Sh(Z)
} \]
est commutatif.

Soit $Z$ et $Z'$ deux espaces sobres munis d'une $G$-action (à gauche). Soit $\alpha:Sh(G;Z)\rightarrow Sh(G;Z')$ un morphisme au-dessus de $B^{sm}_{G}$, i.e. tel que le diagramme
\[ \xymatrix{
Sh(G;Z) \ar[rd]_{f} \ar[r]^{\alpha} &Sh(G;Z') \ar[d]_{f'}   \\
&B^{sm}_{G}
} \]
commute. On obtient un morphisme
\begin{equation}\label{onemap}
Sh(Z)\simeq Sh(G;Z)/_{f^*EG}\longrightarrow Sh(G;Z')/_{f'^*EG}\simeq Sh(Z').
\end{equation}
qui est $G$-équivariant au sens de (\ref{actionGsurSh}). D'après \cite{SGA4} IV. 4.2.3, il existe une unique application continue $a:Z\rightarrow Z'$ tel que $Sh(a)=\alpha$. De plus, $a$ est un homéomorphisme si et seulement si $\alpha$ est une équivalence. Mais le morphisme de topos (\ref{onemap}) est équivariant, i.e. respecte l'action définie par (\ref{actionGsurSh}). D'après ce qui précède, cette action (à droite) de $G$ sur $Sh(Z)$ et $Sh(Z')$ est induite par l'action à gauche de $G$ sur $Z$ et $Z'$. Il suit que $a$ est $G$-équivariante.
\end{f-proof}
\subsection{}La preuve précédente donne aussi le résultat suivant.
\begin{cor}
Soit $Z$ un espace sobre muni d'une action $r:G\times Z\rightarrow Z$. Le morphisme $f:Sh(G;Z)\rightarrow B^{sm}_G$ définit une action $\rho$ à droite de $G$ sur $Sh(Z)\simeq Sh(G;Z)\times_{B^{sm}_G}\underline{Set}$, et on a $\rho=Sh(\check{r})$.
\end{cor}
\begin{f-proof}
On a $$Sh(G;Z)\times_{B^{sm}_G}\underline{Set}\simeq Sh(G;Z)\times_{B^{sm}_G}B^{sm}_G/_{EG}\simeq Sh(G;Z)/_{f^*EG}\simeq Sh(Z)$$
On a de plus $$Sh(Z)\times G\simeq Sh(G;Z)/_{f^*EG}\times_{Set}\underline{Set}/_{G}\simeq Sh(G;Z)/_{f^*EG\times G}$$
et l'action $\rho:Sh(Z)\times G\rightarrow Sh(Z)$ est simplement donnée par le morphisme
$$Sh(G;Z)/_{f^*EG\times G}\longrightarrow Sh(G;Z)/_{f^*EG}$$
lui-même induit par l'action à droite ${f^*EG}\times G\rightarrow{f^*EG}$.
\end{f-proof}
La situation est la même au-dessus d'un topos de base arbitraire. Soit $G$ un groupe topologique. On note encore $G$ le groupe de $\mathcal{T}$ qu'il définit (on pourra plus généralement considérer un topos quelconque $\mathcal{T}$ dont $G$ est un groupe). Soit de plus $\mathcal{S}\rightarrow\mathcal{T}$ un $\mathcal{T}$-topos et $f:\mathcal{S}\rightarrow B_G$ un morphisme au-dessus de $\mathcal{T}$. On note $$\widetilde{\mathcal{S}}:=\mathcal{S}\times_{B_G}\mathcal{T}\simeq \mathcal{S}/_{f^*EG}.$$
Le morphisme $f$ définit une $\mathcal{T}$-action à droite de $G$ sur le $\mathcal{T}$-topos $\widetilde{\mathcal{S}}$. Cette $\mathcal{T}$-action est donnée par le $\mathcal{T}$-morphisme
\begin{equation}\label{action}
\widetilde{\mathcal{S}}\times_{\mathcal{T}}\mathcal{S}(G)\longrightarrow \widetilde{\mathcal{S}},
\end{equation}
et par des isomorphismes de transitivité rendant commutatifs les diagrammes d'une action. On note ici $\mathcal{S}(G):=\mathcal{T}/_G$ et on a
$\widetilde{\mathcal{S}}\times_{\mathcal{T}}\mathcal{S}(G)\simeq \widetilde{\mathcal{S}}/_G$. Le morphisme (\ref{action}) est induit par l'action à droite $EG\times G\rightarrow EG$.

Le $\mathcal{T}$-topos $\widetilde{\mathcal{S}}$ muni de la $\mathcal{T}$-action (\ref{action}) permet de retrouver le topos $\mathcal{S}$ et le morphisme $\mathcal{S}\rightarrow B_G$. En effet, $\mathcal{S}$ est équivalent au topos de descente du topos simplicial
\begin{equation}\label{top-simpl-1}
\widetilde{\mathcal{S}}\times_{\mathcal{T}}\mathcal{S}(G)\times_{\mathcal{T}}\mathcal{S}(G)
\rightrightarrows\rightarrow\widetilde{\mathcal{S}}\times_{\mathcal{T}}\mathcal{S}(G)
\rightrightarrows\leftarrow\widetilde{\mathcal{S}}
\end{equation}
car ce topos simplicial s'identifie à
\begin{equation}\label{top-simpl-2}
\mathcal{S}/_{EG\times EG\times EG}\rightrightarrows\rightarrow\mathcal{S}/_{EG\times EG}
\rightrightarrows\leftarrow\mathcal{S}/_{EG},
\end{equation}
via l'isomorphisme $EG\times EG\simeq EG\times G$ dans $B_G$, où $G$ désigne l'objet de $B_G$ donné par $G$ avec l'action triviale. Le morphisme $\widetilde{\mathcal{S}}\rightarrow\mathcal{T}$ induit un morphisme du topos simplicial (\ref{top-simpl-1})  dans
\begin{equation}\label{top-simpl-3}
\mathcal{T}/_{G\times G}\rightrightarrows\rightarrow\mathcal{T}/_{G}
\rightrightarrows\leftarrow\mathcal{T}
\end{equation}
Le topos de descente de (\ref{top-simpl-3}) est $B_G$ (car $B_G/_{EG}\simeq\mathcal{T}$). Le morphisme du topos simplicial (\ref{top-simpl-1}) dans (\ref{top-simpl-3}) induit un morphisme entre les topos de descente : on retrouve le morphisme $\mathcal{S}\rightarrow B_G$. On peut d'ailleurs montrer le résultat suivant par descente : \emph{la donnée d'un $\mathcal{T}$-topos $\mathcal{S}$ muni d'un morphisme dans $B_G$ équivaut à la donnée du $\mathcal{T}$-topos $\mathcal{S}\times_{B_G}\mathcal{T}$ muni d'une $\mathcal{T}$-action à droite de $G$}.

Dans le cas particulier où $\mathcal{T}$ est le topos des espaces topologiques, $\mathcal{S}(G)\simeq TOP(G)$ est le gros topos de l'espace topologique $G$. Si de plus $\mathcal{S}=\mathcal{S}(G,M)$ où $G$ opère continûment sur un espace $M$, alors la $\mathcal{T}$-action (\ref{action}), déduite du morphisme $\mathcal{S}(G,M)\rightarrow B_G$, est donnée par le morphisme $\mathcal{S}(M\times G)\rightarrow \mathcal{S}(M)$ induit par l'action de $G$ sur $M$.


\begin{thebibliography}{29}

\bibitem{Artin-Verdier}{M. Artin and J.-L. Verdier, \emph{Seminar on Etale Cohomology of Number Fields.} Lecture
notes prepared in connection with a seminar held at the Woods Hole Summer Institute
on Algebraic Geometry, July 6-31, 1964.}
\bibitem{Deninger-some-analogies}{C. Deninger, \emph{Some analogies between number theory and dynamical systems on
foliated spaces.} Doc. Math. J. DMV. Extra Volume ICM I (1998),
23-46}
\bibitem{Deninger-possible-significance}{C. Deninger, \emph{On dynamical systems and their possible significance for arithmetic
geometry.}  Regulators in analysis, geometry and number theory,
29--87, Progr. Math., 171, Birkhäuser Boston, Boston, MA, 2000.}
\bibitem{Deninger-NTandDSonFS}{C. Deninger, \emph{Number theory and dynamical systems on foliated spaces.}  Jahresber.
Deutsch. Math.-Verein.  103  (2001), no. 3, 79--100.}
\bibitem{Deninger2}{C. Deninger, \emph{A note on arithmetic topology and dynamical systems.} Algebraic number theory and algebraic geometry, 99-114,
Contemp. Math., 300, Amer. Math. Soc., Providence, RI, 2002.}
\bibitem{Deninger-explicit formulas}{C. Deninger, \emph{On the nature of the "explicit formulas" in
analytic number theory, a simple example}.  Number theoretic methods
(Iizuka, 2001),  97--118, Dev. Math., 8, Kluwer Acad. Publ.,
Dordrecht, 2002.}
\bibitem{Deninger3}{C. Deninger, \emph{Analogies between analysis on foliated spaces and arithmetic geometry.}  Groups and analysis,  174--190, London Math. Soc. Lecture Note Ser., 354, Cambridge Univ. Press, Cambridge, 2008.}
\bibitem{Deninger4}{C. Deninger, \emph{The Hilbert-Polya strategy and height pairing.} Preprint, 2010.}
\bibitem{Flach-moi}{M. Flach; B. Morin, \emph{On the Weil-\'etale topos of regular arithmetic schemes.} Preprint (2010).}
\bibitem{SGA4}{A. Grothendieck, M. Artin and J.L. Verdier, \emph{Théorie des Topos et cohomologie étale des schémas
(SGA4).} Lectures Notes in Math. 269, 270, 305, Springer, 1972.}
\bibitem{Lichtenbaum-finite-field}{S. Lichtenbaum, \emph{The Weil-étale topology on schemes over finite fields.}
Compositio Math.  141  (2005),  no. 3, 689--702.}
\bibitem{Lichtenbaum}{S. Lichtenbaum, \emph{The Weil-étale topology for Number Rings.}  Ann. of Math (2) 170 (2009),  no. 2, 657-683.}
\bibitem{Leichtnam-invitation}{E. Leichtnam, \emph{An invitation to Deninger's work on arithmetic zeta functions.}  Geometry, spectral theory, groups, and dynamics,  201--236, Contemp. Math., 387, Amer. Math. Soc., Providence, RI, 2005.}
\bibitem{Leichtnam}{E. Leichtnam, \emph{Scaling group flow and Lefschetz trace formula for laminated spaces with $p$-adic transversal.} Bull. Sci. Math.  131  (2007),  no. 7, 638--669.}
\bibitem{Milne-motivesfinitefields}{J. Milne, \emph{Motives over finite fields.}  Motives (Seattle, WA, 1991),  401--459, Proc. Sympos. Pure Math., 55, Part 1.}
\bibitem{moi}{B. Morin, \emph{Utilisation d'une cohomologie étale équivariante en topologie arithmétique.} Compositio Math. 144
(2008), no. 1, 32-60.}
\bibitem{these}{B. Morin, \emph{Sur le topos Weil-étale d'un corps de nombres.} Thesis (2008).}
\bibitem{Fund-group-I}{B. Morin, \emph{The Weil-étale fundamental group of a number field I.} Preprint (2009).}
\bibitem{Morishita-primesandknots}{M. Morishita, \emph{Analogies between prime numbers and knots.}
Sugaku 58 (2006), no. 1, 40-63.}
\end{thebibliography}
\end{document}